\documentclass[11pt]{amsart}
\pdfoutput=1

\usepackage{mathptmx}
\usepackage{amssymb, amsfonts, amsthm}
\usepackage{tikz}
\usepackage{pgf}
\usepackage{epstopdf}
\usepackage{booktabs}
\usepackage{multicol}
\usepackage{placeins}
\usepackage{caption} 
\usepackage{graphicx}
\usepackage{mathdots}
\usepackage[capitalise]{cleveref}
\usepackage{mathrsfs}
\usepackage{txfonts}
\usepackage{amsopn}
\usepackage{natbib}
\usepackage{graphicx}
\usepackage{amsmath,amsthm,bm,mathrsfs}

 \newtheoremstyle{theorem}{6pt}{6pt}{\rm}{}{\sffamily}{ }{ }{}
 \theoremstyle{theorem}
\newtheorem{theorem}{\sc Theorem}[section]

 \newtheoremstyle{algorithm}{6pt}{6pt}{\rm}{}{\sffamily}{ }{ }{}
 \theoremstyle{algorithm}
\newtheorem{algorithm}{\sc Algorithm}[section]

 \newtheoremstyle{lemma}{6pt}{6pt}{\rm}{}{\sffamily}{ }{ }{}
 \theoremstyle{lemma}
 \newtheorem{lemma}{\sc Lemma}[section]

\newtheoremstyle{case}{6pt}{6pt}{\rm}{}{\sffamily}{. }{ }{}
 \theoremstyle{case}

 \newtheoremstyle{statement}{6pt}{6pt}{\rm}{}{\sffamily}{ }{ }{}
\theoremstyle{statement}

 \newtheoremstyle{corollary}{6pt}{6pt}{\rm}{}{\sffamily}{ }{ }{}
 \theoremstyle{corollary}
 \newtheorem{corollary}{\sc Corollary}[section]
 
  \newtheoremstyle{definition}{6pt}{6pt}{\rm}{}{\sffamily}{ }{ }{}
 \theoremstyle{definition}
 \newtheorem{definition}{\sc Definition}[section]

\newtheoremstyle{example}{6pt}{6pt}{\rm}{}{\sffamily}{ }{ }{}
\theoremstyle{example}

\newtheoremstyle{remark}{6pt}{6pt}{\rm}{}{\sffamily}{ }{ }{}
\theoremstyle{remark}
\newtheorem{remark}{\sc Remark}[section]

\newtheoremstyle{approximation}{6pt}{6pt}{\rm}{}{\sffamily}{ }{ }{}
\theoremstyle{approximation}

\newtheoremstyle{scheme}{6pt}{6pt}{\rm}{}{\sffamily}{ }{ }{}
\theoremstyle{scheme}

\newtheoremstyle{Algorithm}{6pt}{6pt}{\rm}{}{\sffamily}{ }{ }{}
\theoremstyle{Algorithm}

\newtheoremstyle{Assumption}{6pt}{6pt}{\rm}{}{\sffamily}{ }{ }{}
\theoremstyle{Assumption}

\newtheoremstyle{proposition}{6pt}{6pt}{\rm}{}{\sffamily}{ }{ }{}
\theoremstyle{proposition}

\newtheoremstyle{hypo}{6pt}{6pt}{\rm}{}{\sffamily}{ }{ }{}
 \theoremstyle{hypo}

  \newtheoremstyle{Step}{6pt}{6pt}{\rm}{}{}{ }{ }{}
 \theoremstyle{Step}

\numberwithin{equation}{section}

\usepackage{amsfonts,amssymb}
\usepackage{subfigure,epsfig,psfrag}
\usepackage{epstopdf}
\usepackage{color,url}
\usepackage{algorithm,algorithmic}
\usepackage{tikz}
\usepackage{graphicx}

\newcommand{\mK}{\mathcal{K}}
\newcommand{\mI}{\mathcal{I}}
\newcommand{\mJ}{\mathcal{J}}

\newcommand{\mM}{\mathcal{M}}

\newcommand{\RR}{\mathbb{R}}

\newcommand{\NN}{\mathbb{N}}
\newcommand{\CC}{\mathbb{C}}

\def\Dpartial#1#2{ {\partial #1 \over \partial #2} }

\newcommand\norm[1]{\left\Vert{} #1 \right\Vert}

\def\J{\mathcal{J}}

\newcommand{\argmin}{\operatorname{argmin}}

\begin{document}

\title[Breast Cancer Detection through EIT and Optimal Control Theory]{Breast Cancer Detection through Electrical Impedance Tomography and Optimal Control Theory:
        Theoretical and Computational Analysis}
	 \thanks{Department of Mathematics, Florida Institute of Technology, Melbourne, FL 32901}

	\author[U. G. Abdulla]{Ugur G. Abdulla}
	\address{Department of Mathematics, Florida Institute of Technology, Melbourne, FL 32901}
	\email{abdulla@fit.edu}	
	
	\author[V. Bukshtynov]{Vladislav Bukshtynov}	
    %       \address{Department of Mathematics, Florida Institute of Technology, Melbourne, FL 32901}
	%\email{vbukshtyn@fit.edu}	
		
	\author[S. Seif]{Saleheh Seif}

\maketitle

% abstract style
% Two grouping braces are necessary in abstract environment:
% first argument contains abstract text; second argument contains keywords text

\begin{abstract}
{The Inverse Electrical Impedance Tomography (EIT) problem on recovering electrical conductivity tensor and potential in the body based on the measurement of the boundary voltages on the electrodes for a given electrode current is analyzed. A PDE constrained optimal control framework in Besov space is pursued, where the electrical conductivity tensor and boundary voltages are control parameters, and the cost functional is the norm declinations of the boundary electrode current from the given current pattern and boundary electrode voltages from the measurements. The state vector is a solution of the second order elliptic PDE in divergence form with bounded measurable coefficients under mixed Neumann/Robin type boundary condition. Existence of the optimal control and Fr\'echet differentiability in the Besov space setting is proved. The formula for the Fr\'echet gradient and optimality condition is derived. Extensive numerical analysis is pursued in the 2D case by implementing the projective gradient method, re-parameterization via principal component analysis (PCA) and Tikhonov regularization. 
}
{Breast cancer detection, Electrical Impedance Tomography, PDE constrained optimal control, Fr\'echet differentiability, projective gradient method, principal component analysis, Tikhonov regularization.}
\end{abstract}

\section{Introduction and Problem Description}
\label{sec:Introduction}

This paper analyzes mathematical model for the breast cancer detection through EIT.  Let $Q \in \mathbb{R}^n$ be an open and bounded set  representing body, and assume $ A(x) = \big( a_{ij}(x) \big)_{ij=1}^n$ be a matrix representing the electrical conductivity tensor at the point $x\in Q$. Electrodes, $(E_l)_{l=1}^m$, with contact impedances vector $Z:= (Z_l)_{l=1}^m \in \mathbb{R}^m_+$ are attached to the periphery of the body, $\partial Q$. Electrical currents vector $I:=(I_l)_{l=1}^m\in \mathbb{R}^m$ is applied to the electrodes. Vector $I$ is called {\it current pattern} if it satisfies conservation of charge condition
\begin{equation}\label{currentpattern}
\sum_{l=1}^m I_l=0
\end{equation}
The induced constant voltage on electrodes is denoted by $ U:= (U_l)_{l=1}^m \in \mathbb{R}^m$. By specifying ground or zero potential it is assumed that
\begin{equation}\label{grounding}
\sum_{l=1}^m U_l=0
\end{equation}
EIT problem is to find the electrostatic potential $u: Q\to \mathbb{R}$ and boundary voltages $ U$ on  $(E_l)_{l=1}^m $. The mathematical model of the EIT problem is described through the following boundary value problem for the second order elliptic partial differential equation:

\begin{alignat}{3}
& -\sum_{i,j=1}^n \big( a_{ij}(x) u_{x_j}\big)_{x_i} = 0, & \quad  x\in Q    \label{ general pde} \\
& \frac{\partial u(x)}{\partial \mathcal{N}}  = 0,       & \quad  x \in \partial Q -\bigcup\limits_{l=1}^{m} E_{l}\label{neumann}\\
& u(x) +  Z_l \frac{\partial  u(x)}{\partial \mathcal{N}}  = U_{l} ,       &\quad  x\in E_{l}, \ l= \overline{1,m}\label{bdry cond. of gen pde}\\
& \displaystyle\int_{E_l}  \frac{\partial  u(x)} {\partial \mathcal{N}} ds = I_l,       &\quad  l= \overline{1,m} \label{boundaryflux}
\end{alignat}
where
$$\frac{\partial u(x)}{\partial \mathcal{N}} = \sum_{i,j}a_{ij}(x)u_{x_j} \nu^i$$
be a co-normal derivative at $x$, and $\nu=(\nu^1,...,\nu^n)$ is the outward normal at a point $x$ to $\partial Q$.
Electrical conductivity matrix $A=(a_{ij})$ is positive definite with
\begin{equation}
\displaystyle\sum_{i,j=1}^n a_{ij}(x) \xi_i \xi _j  \geq \mu \displaystyle\sum_{i=1}^n \xi_i^2 , \ \forall \xi\in \mathbb{R}^n;
\quad \mu >0.
\label{reconstructive. inequality.}
\end{equation}
%$\nu=(\nu^1,...,\nu^n)$ is the outward normal at a point $x$ to $\partial Q$,
%\[  \mathcal{N} = \Big (\sum_{i}a_{i1}(x) \nu^i,...,\sum_{i}a_{in}(x) \nu^i\Big ) \]
%s a co-normal vector at the point $x$ and
The following is the

{\bf EIT Problem:} {\it Given electrical conductivity tensor $A$, electrode contact impedance vector $Z$, and electrode {\it current pattern} $I$ it is required to find electrostatic potential $u$ and electrode voltages $U$ satisfying \eqref{grounding}--\eqref{boundaryflux}}:
\[ (A,Z,I) \longrightarrow (u,U) \]
The goal of the paper is to analyze inverse EIT problem of determining conductivity tensor $A$ from the measurements of the boundary voltages $U^*$.

{\bf Inverse EIT Problem:} {\it Given electrode contact impedance vector $Z$, electrode {\it current pattern} $I$ and boundary electrode measurement $U^*$, it is required to find electrostatic potential $u$ and electrical conductivity tensor $A$ satisfying \eqref{grounding}--\eqref{boundaryflux}} with $U=U^*$.

Mathematical model  \eqref{grounding}--\eqref{boundaryflux} for the EIT Problem was suggested in \cite{cheng1989electrode}. In \cite{somersalo1992existence} it was demonstrated that the model is capable of predicting the experimentally measured voltages to within 0.1 percent. Existence and uniqueness of the solution to the problem  \eqref{grounding}-\eqref{boundaryflux} was also proved in \cite{somersalo1992existence}. EIT is a rapidly developing non-invasive imaging technique recently gaining popularity in various medical applications including breast screening and cancer detection \cite{zou2003review,brown2003electrical,adlergreit,holder2004electrical}. The objective of the Inverse EIT Problem is reconstructing the electrical conductivity through measuring voltages of electrodes placed on the surface of a test volume. The electrical conductivity of the malignant tumors of the breast may significantly differ from the conductivity of surrounding normal tissue. This provides a possible way to develop an efficient, safe and inexpensive method to detect and localize such tumors.
X-ray mammography, ultrasound, and magnetic resonance imaging (MRI) are among methods that are used currently for breast cancer diagnosis \cite{zou2003review}. However, these methods have various flaws and cannot distinguish breast cancer from benign breast lesions with certainty \cite{zou2003review}. EIT is a fast, inexpensive, portable, and relatively harmless technique, although it also has the disadvantage of poor image resolution \cite{paulson1995pompus}.  Different types of regularization have been applied to overcome this issue \cite{brown2003electrical,adler2005eidors}.
Inverse EIT Problem is an ill-posed problem and belongs to the class of so-called Calderon type inverse problems, due to pioneering work \cite{calderon} where well-posedness of the inverse problem for the identification of the conductivity coefficient of the second order elliptic PDE through Dirichlet-to-Neumann or Neumann-to-Dirichlet maps is presented. We refer to topical review
paper \cite{borcea} on EIT and Calderon type inverse problems. Reconstruction of the coefficient in Calderon problem is pursued in \cite{nachman1988reconstructions},and the uniqueness of the solution has been demonstrated \cite{sylvester1987global}. This framework was shown to be stable in \cite{alessandrini1988stable}. Well-posedness of the inverse Calderon problem with partial boundary data is analyzed in \cite{kenig2007}.
Statistical methods have been applied for solving inverse EIT problem in \cite{kaipio2000statistical,kaipio1999inverse,roininen2014whittle}. Bayesian formulation of EIT in infinite dimensions has been proposed in \cite{dunlop2015bayesian}.
An experimental iterative algorithm, POMPUS, was introduced, the accuracy of which is comparable to standard Newton-based algorithms \cite{paulson1995pompus}.
An analytic solution for potential distribution on a 2D homogeneous disk for EIT problem was analyzed in \cite{demidenko2011analytic}. A statistical model called gapZ, has also been developed for solving EIT using Toeplitz matrices \cite{demidenko2011statistical}.

In this paper, inverse EIT Problem is investigated with unknown electrical conductivity tensor $A$. This is in contrast with current state of the art in the field where usually inverse EIT problem is solved for the reconstruction of the single conductivity function. This novelty is essential in understanding and detection of the highly anisotropic distribution of the cancerous tumor in breast. We formulate Inverse EIT Problem as a PDE constrained optimal control problem in Besov spaces framework, where the electrical conductivity tensor and boundary voltages are control parameters, and the cost functional is the norm declinations of the boundary electrode current from the given current pattern and boundary electrode voltages from the measurements.
We prove the existence of the optimal control and Fr\'echet differentiability in the Besov space setting. The formula for the Fr\'echet gradient and optimality condition is derived. Based on the Fr\'echet differentiability result we develop projective gradient method in Besov spaces. Extensive numerical analysis in the 2D case by implementing the projective gradient method, re-parameterization via PCA and Tikhonov regularization is pursued.

The organization of the paper is as follows. In Section~\ref{sec:notations} we introduce the notations of the functional spaces. In Section~\ref{sec:Optimal_Control_Problem} we introduce Inverse EIT Problem as PDE constrained optimal control problem. In Section~\ref{sec:Main_Results} we formulate the main results.
Proof of the main results are presented in Section~\ref{sec:Proofs_Main_Results}. In Section~\ref{sec:numerical_results} we present the results of the computational analysis for the 2D model. Finally, in Section~\ref{sec:conclusions} we outline the main conclusions.

\section{Notations}
\label{sec:notations}
In this section, assume $Q$ is a domain in $\mathbb{R}^n$.
\begin{itemize}

\item For $1 \leq p < \infty$,  $L_p(Q)$ is a Banach space of measurable functions on $Q$ with finite norm

$$
\| u\|_{L_p(Q)} := \Big( \int_{Q} |u(x)|^p dx \Big)^{\frac{1}{p}}
$$
In particular if $p=2$, $L_2(Q)$ is a Hilbert space with inner product 

\begin{equation*}
(f,g)_{L_2(Q)} = \int_{Q} f(x)g(x) dx
\end{equation*}
\item $L_{\infty}(Q)$ is a Banach space of measurable functions on $Q$ with finite norm 
\[
\|u\|_{L_{\infty}(Q)} := \text{ess} \displaystyle\sup_{x \in Q} |u(x)|
\]

\item For $s \in \mathbb{Z}_+$,  $W_p^s(Q)$ is the Banach space of measurable functions on $Q$ with finite norm
$$
\| u\|_{W_p^s(Q)} := \Big(\int_Q \displaystyle\sum_{|\alpha| \leq s}  |D^{\alpha}u(x)|^p dx \Big)^{\frac{1}{p}},
$$
where $\alpha = (\alpha_1,\alpha_2,..., \alpha_n )$,  $\alpha_j$ are nonnegative integers, $|\alpha| = \alpha_1+...+\alpha_n$, $D_k = \frac{\partial }{\partial x_k}$,  $ D^{\alpha}= D_1^{\alpha_1}...D_n^{\alpha_n}. $
In particular if $p=2$, $H^s(Q) := W_2^s(Q)$  is a Hilbert space with inner product
\begin{equation*}
(f,g)_{H^s(Q)} = \displaystyle\sum_{|\alpha| \leq s}  (D^{\alpha}f(x),D^{\alpha}g(x) )_{L_2(Q)}
\end{equation*}

\item For $s \notin \mathbb{Z}_+$, $B_p^s(Q)$ is the Banach space of measurable functions on $Q$ with finite norm

\begin{equation*}
\|u\|_{B_p^s(Q)} := \| u\|_{W_p^{[s]}(Q)} + [u]_{B_p^s(Q)}
\end{equation*}

where 
\begin{equation*}
 [u]_{B_p^s(Q)} := \int_Q \int_Q \frac{\big|\frac{\partial^{[s]} u(x)}{\partial x^{[s]} } -\frac{\partial^{[s]} u(y)}{\partial x^{[s]} }\big|^p}{|x-y|^{1+p(s - [s])}} dx dy\Big)^{\frac{1}{p}}
\end{equation*}

$H^{\epsilon}(Q) := B_2^{\epsilon}(Q)$ is an Hilbert space.
\item $\textbf{ba} (Q) = \big(L_{\infty}(Q)\big)'$ is the  Banach space of bounded and finitely additive signed measures on $Q$ and the dual space of $L_{\infty}(Q)$ with finite norm

$$\|\phi \|_{\textbf{ba} (Q)} = | \phi |(Q),$$

$|\phi|(Q) $ is total variation of $\phi$ and defined as $ | \phi |(Q) = \sup \displaystyle\sum_i \phi (E_i) $, where the supremum is taken over all partitions  $\cup E_i$ of $E$ into measurable subsets $E_i$.
\item $\mathbb{M}^{m\times n}$ is a space of real $m\times n$ matrices.
\item $\mathscr{L} :=L_{\infty} (Q; \mathbb{M}^{n \times n})$ is the Banach space of $n \times n$ matrices of $L_{\infty}(Q)$ functions.
\item $ \mathscr{L}':=  \textbf{ba} (Q; \mathbb{M}^{n \times n}) = \bigl( L_{\infty} (Q; \mathbb{M}^{n \times n})  \bigr)' $ is the Banach space of $n \times n$ matrices of $\textbf{ba} (Q)$ measures.
\end{itemize}
\section{Optimal Control Problem}
\label{sec:Optimal_Control_Problem}
We formulate Inverse EIT Problem as the following PDE constrained optimal control problem. Consider the minimization of the cost functional
\begin{align}
  \mathcal{J}(v) =   \displaystyle\sum _{l=1}^m \Big |  \displaystyle\int_{E_l} \frac{U_l-u(x)}{Z_l}ds- I_l\Big|^2 + \beta |U-U^*|^2\label{eq:cost_functional}
\end{align}
on the control set
%%%%%%%%%%%%%%%%%%%
%\begin{gather}
%V_R=\big\{v = (A, U)\in  L_{\infty} (Q; \mathbb{R}^{n \times n})
%\times \mathbb{R}^m \Big| \|A\|_{L_{\infty}}+ |U|_{\mathbb{R}^m} \leq R
%\nonumber\\
%\sum_{l=1}^m U_l = 0, \  \xi^T A \xi \geq \nu |\xi|^2, \  \forall \xi\in \mathbb{R}^n, \  \nu \geq 0   \big \}\nonumber
%\end{gather}
%%%%%%%%%%%%%%%%%%%%%%
\begin{gather}
V_R=\big\{v = (A, U)\in \Big( L_{\infty} (Q; \mathbb{M}^{n \times n}) \bigcap H^{\epsilon}(Q; \mathbb{M}^{n \times n}) \Big )
\times \mathbb{R}^m \Big| \sum_{l=1}^m U_l = 0,
\nonumber\\
\|A\|_{L_{\infty}}+ \|A\|_{H^{\epsilon}} +|U| \leq R, \  \xi^T A \xi \geq \mu |\xi|^2, \  \forall \xi\in \mathbb{R}^n, \  \mu > 0  \big \}\nonumber
\end{gather}
where $\beta> 0$, and $u=u(\cdot;v)\in H^1(Q)$ is a solution of the elliptic problem \eqref{ general pde}--\eqref{bdry cond. of gen pde}.
This optimal control problem will be called Problem $\mJ$. The first term in the cost functional $\mathcal{J}(v)$ characterizes
the mismatch of the condition \eqref{boundaryflux} in light of the Robin condition \eqref{bdry cond. of gen pde}.

Note that the variational formulation of the EIT Problem is a particular case of the Problem $\mJ$, when the conductivity tensor $A$ is known, and therefore is removed from the control set by setting $R=+\infty$ and $\beta=0$:
\begin{align}
  \mathcal{I}(U) =   \displaystyle\sum _{l=1}^m \Big |  \displaystyle\int_{E_l} \frac{U_l-u(x)}{Z_l}ds- I_l\Big|^2\label{eq:cost_functional1}\to \inf
\end{align}
in a control set
\begin{equation}
W=\big\{U\in \mathbb{R}^m \Big| \sum_{l=1}^m U_l = 0 \big\}
\end{equation}
where $u=u(\cdot;v)\in H^1(Q)$ is a solution of the elliptic problem \eqref{ general pde}--\eqref{bdry cond. of gen pde}.
This optimal control problem will be called Problem $\mI$. It is a convex PDE constrained optimal control problem (Remark~\ref{convex}, Section~\ref{sec:Proofs_Main_Results}). 

Inverse EIT problem on the identification of the electrical conductivity tensor $A$ with $m$ input data $(I_l)_{l=1}^m$ is highly ill-posed.
Next, we formulate an optimal control problem which is adapted to the situation when the size of input data can be increased 
through additional measurements while keeping the size of the unknown parameters fixed. Let $U^1=U, I^1=I$ and consider
$m-1$  new permutations of boundary voltages
\begin{equation}\label{permutations}
U^j=(U_j,...,U_{m},U_1,...,U_{j-1}), \ j=2,...,m
\end{equation}
applied to electrodes $E_1, E_2, \ldots, E_{m}$ respectively. Assume that the ``voltage--to--current'' measurement allows us 
to measure associated currents $I^j=(I^j_1,\cdots,I^j_{m})$. By setting $U^1=U^*$ and having a new set of $m^2$ input data $(I^j)^{m}_{j=1}$,
we now consider optimal control problem on the minimization of the new cost functional
\begin{align}
  \mathcal{K}(v) = \displaystyle\sum_{j=1}^m \displaystyle\sum_{l=1}^m
  \Big | \displaystyle\int_{E_l} \frac{U^j_l-u^j(x)}{Z_l} ds - I^j_l\Big|^2+\beta|U-U^*|^2
  \label{eq:cost_functional_rotation}
\end{align}
on a control set $V_R$, where each function $u^j(\cdot; A, U^j), j=1,...,m$, solves elliptic PDE problem \eqref{eq:forward_1}--\eqref{eq:forward_3}
with $U$ replaced by $U^j$. This optimal control problem will be called Problem $\mK$. %It should be noted that similar approach can
%be pursued to increase the size of input data to $m!$ by adding all possible permutations of $U$ in \eqref{permutations}.

We effectively use Problem $\mI$ to generate model examples of the inverse EIT problem which adequately represents the diagnosis
of the breast cancer in reality. Computational analysis based on the Fr\'echet differentiability result and gradient method in Besov spaces for the Problems $\mJ$ and $\mK$ is pursued in realistic model examples. 
\section{Main Results}
\label{sec:Main_Results}
Let bilinear form $B: H^1(Q)\times H^1(Q)\to \mathbb{R}$ be defined as
\begin{equation}
B[u, \eta] = \int_Q  \displaystyle\sum_{i,j=1}^n a_{ij}u_{x_j} \eta_{x_i}  dx + \displaystyle\sum_{l=1}^m \frac{1}{Z_l}\int_{E_l} u \eta ds,
\end{equation}
\begin{definition}\label{weaksolution}
For a given $v\in V_R$, $u=u(\cdot;v)\in H^1(Q)$ is called a solution of the problem \eqref{ general pde}--\eqref{bdry cond. of gen pde} if
\begin{equation}\label{weaksolution}
B[u, \eta] = \displaystyle\sum_{l=1}^m \frac{1}{Z_l}\int_{E_l} \eta  U_l ds, \quad  \forall \eta \in H^1(Q).
\end{equation}
\end{definition}
For a given control vector $v\in V_R$ and corresponding $u(\cdot;v)\in H_1(Q) $, consider the adjoined problem:
\begin{alignat}{3}
& \displaystyle\sum_{ij}   (a_{ij} \psi_{x_i})_{x_j}  = 0, \quad & x \in Q
\label{adj pde}
\\
& \frac{\partial \psi}{ \partial \mathcal{N}} = 0 , & x \in \partial Q -\bigcup\limits_{l=1}^{m} E_{l}
\\
& \psi + Z_l \frac{\partial \psi}{ \partial \mathcal{N}}=  2\int_{E_l} \frac{u-U_l}{Z_l} ds + 2 I_l ,       \quad  & x\in E_{l}, \ l= \overline{1,m}
\label{bdry cond of adj pde}
\end{alignat}
\begin{definition}\label{weakadjoinedsolution}
$\psi \in H^1(Q)$ is called a solution of the adjoined problem \eqref{adj pde}--\eqref{bdry cond of adj pde} if
\begin{equation}\label{weaksolutionadjoined}
B[\psi, \eta] = \displaystyle\sum_{l} \int_{E_l} \frac{\eta}{Z_l} \Big[  2\int_{E_l} \frac{u-U_l}{Z_l} ds + 2 I_l \Big] ds,\  \forall \eta \in H^1(Q).
\end{equation}
\end{definition}
In Lemma~\ref{existence lax}, Section~\ref{sec:Proofs_Main_Results} it is demonstrated that for a given $v\in V_R$, both elliptic problems are uniquely solvable.
\begin{definition}\label{defn:Frechet}
    Let $V$ be a convex and closed subset of the Banach space $H$.
    We say that the functional $\mathcal{J}:V\to \mathbb{R}$ is differentiable in the sense of Fr\'echet at the point $v\in V$ if there exists an element $\mathcal{J}'(v) \in H'$ of the dual space such that
    \begin{equation}\label{frechetderivativedefinition}
        \mathcal{J}(v+h)-\mathcal{J}(v)=\left\langle{}\J'(v),{h}\right\rangle_H+o(h,v),
    \end{equation}
    where $v+h\in V\cap\{u: \norm{u}< \gamma\}$ for some $\gamma>0$; $\langle{} \cdot,\cdot\rangle{}_H$ is a pairing
    between $H$ and its dual $H'$, and
    \[
        \frac{o(h,v)}{\norm{h}} \to 0,\quad \text{as}~\norm{h} \to 0.
    \]
    The expression $d\mathcal{J}(v)=\left\langle{}\mathcal{J}'(v),{\cdot}\right\rangle_H$ is called a Fr\'echet differential of $\mathcal{J}$ at $v\in V$, and the element $\mathcal{J}'(v) \in H'$ is called Fr\'echet derivative or gradient of $\J$ at $v\in V$.
\end{definition}
Note that if Fr\'echet gradient $\mathcal{J}'(v)$ exists at $v\in V$, then the Fr\'echet differential $d\mathcal{J}(v)$ is uniquely defined on a convex cone (\cite{abdulla13, abdulla15, abdulla2017frechet, abdulla2018frechet, abdulla2019cam})
\[
    \mathcal{H}_v
    = \{w\in H: w=\lambda (u-v), \lambda \in [0,+\infty), u \in V\}.
\]
The following are the main theoretical results of the paper:
\begin{theorem}(Existence of an Optimal Control).  \label{thm:continuity_of_J}
Problem ${\mathcal J}$ has a solution, i.e.
\begin{equation}\label{solutionexists}
V_* = \{  v = (A, U) \in  V_R; \mathcal{J}(v) = \mathcal{J}_* = \inf_{v \in V_R} \mathcal{J}(v)  \} \neq \emptyset
\end{equation}
\end{theorem}
\begin{theorem}(Fr\'echet Differentiability): \label{thm:Frechet_Differentiability}
The functional $\mathcal{J}(v)$ is differentiable on $V_R$ in the sense of Fr\'echet; the Fr\'echet differential $d\mathcal{J}(v)$ and the gradient $\mathcal{J}'( A , U)\in \mathscr{L}' \times \mathbb{R}^m$ are
\begin{gather}
\left\langle{}\mathcal{J}'(v),\delta v\right\rangle_H=-\int_Q\sum_{i,j=1}^{n}u_{x_j}\psi_{x_i}\delta a_{ij}dx\nonumber\\
+ \sum_{k=1}^m\Big (\displaystyle\sum_{l=1}^m 2 \Big [\int_{E_l} \frac{U_l - u}{Z_l} ds - I_l\Big ]  \int_{E_l} \frac{1}{Z_l} (\delta_{lk} - w^k(s))ds + 2 \beta (U_k - U_k^*)\Big )\delta U_k\\
\mathcal{J}'( A , U ) = \Big(\mathcal{J}'_A ( A , U ), \mathcal{J}'_U( A , U ) \Big)
\nonumber\\
=\left( -  \bigr(\psi_{x_i}  u_{x_j} \bigl)_{i,j=1}^n ,  \Big( \displaystyle\sum_{l=1}^m 2 \Big[\int_{E_l} \frac{U_l - u}{Z_l} ds - I_l\Big]  \int_{E_l} \frac{1}{Z_l} (\delta_{lk} - w^k(s))ds + 2 \beta (U_k - U_k^*)  \Big)_{k=1}^m  \right) \label{eq:Frechet_Derivative}
\end{gather}
where $u=u(\cdot;v), \psi=\psi(\cdot;v)$; $w^k = \frac{\partial u}{\partial U_k}=u(\cdot;A,e_k), \ k = 1,2,..,m$ is a solution of \eqref{ general pde}--\eqref{bdry cond. of gen pde} with $v=(A,e_k) $, $e_k\in\mathbb{R}$ is a unit ort vector in $x_k$-direction; $\delta_{lk}$ is a Kronecker delta; $\delta v = (\delta A,\delta U)=((\delta a_{ij})_{i,j=1}^n,(\delta U_k)_{k=1}^m)$ is a variation of the control vector $v \in V_R$ such that $v + {\delta v} \in V_R$.
\end{theorem}
\begin{corollary}(Optimality Condition)
 If $\mathbf{v}\in V_R$ is an optimal control in Problem $\mathcal{J}$,
    then the following variational inequality is satisfied:
\begin{equation}
\big <\mathcal{J}'(\mathbf{v}), v- \mathbf{v} \big >_H \geq 0,\ \forall v \in V_R.
\end{equation}
\label{optimalitycondition}
\end{corollary}
%%%%%%%%%%%%
\begin{corollary}(Fr\'echet Differentiability): \label{cor:Frechet_Diffrentiability_K}
The functional $\mathcal{K}(v)$ is differentiable on $V_R$ in the sense of Fr\'echet and the Fr\'echet gradient $\mathcal{K}'( \sigma , U)\in \mathscr{L}' \times \mathbb{R}^m$ is
\begin{gather}
\mathcal{K}'(v) = \Big(\mathcal{K}'_A (A, U), \mathcal{K}'_U(A, U) \Big)=
\nonumber\\
\left(-\Big(\displaystyle\sum_{j=1}^m\psi^j_{x_p}  u^j_{x_q} \Big)_{p,q=1}^n,  \Big( \displaystyle\sum_{j=1}^m \displaystyle\sum_{l=1}^m 2 \Big[\int_{E_l} \frac{U^j_l - u_j}{Z_l} ds - I^j_l\Big]  \int_{E_l} \frac{\delta_{l,\theta_{kj}} - w^{\theta_{kj}}(s)}{Z_l}ds + 2 \beta (U_k - U_k^*)  \Big)_{k=1}^m  \right) \label{eq:Frechet_Derivative1}
\end{gather}
where $\psi^j(\cdot), j=1,...,m$, be a solution of the adjoined PDE problem \eqref{adj pde}--\eqref{bdry cond of adj pde}
with $u(\cdot), U$ and $I$ replaced with $u^j(\cdot), U^j$, $I^j$ respectively, and
\begin{equation*}
  \theta_{kj}= \left\{
   \begin{aligned}
   & k-j+1, \qquad&\text{if} \ j\leq k,\\
   & m+k-j+1, \qquad&\text{if} \ j>k.
    \end{aligned} \right.
\end{equation*}
\end{corollary}
%%%%%%%%%%%%%%%%%%%%
\subsection{Gradient Method in Banach Space}
\label{sec:Gradient Method in Banach Space}

Fr\'echet differentiability result of Theorem~\ref{thm:Frechet_Differentiability} and the formula \eqref{eq:Frechet_Derivative} for the Fr\'echet derivative suggest
the following algorithm based on the projective gradient method in Banach space $H$ for the Problem $\mathcal{J}$.

\begin{description}
\item[{ Step 1.}] Set $N=0$ and choose initial vector function $(A^0, U^0) \in V_R$ where
\[ A^0 = (a_{ij}^0)_{ij=1}^n, \ U^0=(U_1^0,...,U_m^0), \ \sum\limits_{l=0}^m U^0_l=0 \]
\item[{Step 2.} ] Solve the PDE problem~\eqref{ general pde}--\eqref{bdry cond. of gen pde} to find $u^N=u(\cdot;A^N,U^N)$ and $\mathcal{J}(A^N,U^N)$.
\item[{Step 3.}] If $N=0$, move to Step 4.
  Otherwise, check the following criteria:
  \begin{equation}
    \left| \frac{\mathcal{J}(A^N,U^N)-\mathcal{J}(A^{N-1},U^{N-1})}{\mathcal{J}(A^{N-1},U^{N-1})} \right| <\epsilon, \quad
   \frac{ \norm{A^{N}-A^{N-1}}}{\norm{A^{N-1}}} < \epsilon,\  \frac{ |U^{N}-U^{N-1}|}{|U^{N-1}|} < \epsilon \label{convergencecriteria}
  \end{equation}
  where $\epsilon$ is the required accuracy.
  If the criteria are satisfied, then terminate the iteration. Otherwise, move to Step 4.

  \item[{ Step 4.}] Solve the PDE problem~\eqref{ general pde}--\eqref{bdry cond. of gen pde}  to find $w_k^N=u(\cdot;A^N,e_k), k=1,...,m$,
 \item[{ Step 5.}] Solve the adjoined PDE problem~\eqref{adj pde}--\eqref{bdry cond of adj pde} to find $\psi_N=\psi(\cdot ;A^N,U^N,u^N)$.
\item[{ Step 6.}] Choose stepsize parameter $\alpha_N>0$ and compute a new control vector components
  $\tilde{A}^{N+1} = (\tilde{a}^{N+1}_{ij}(x))_{i,j=1}^n, \tilde{U}^{N+1}\in\mathbb{R}^m$ as follows:
\begin{gather}
  \tilde{a}^{N+1}_{ij}(x)=a^{N}_{ij}(x)+\alpha_N \psi^N_{x_i}  u^N_{x_j}, \ i,j=1,...,n,\label{gradientupdate_A}
  \\
  \tilde{U}^{N+1}_k=  U^N_k-\alpha_N\Big [ \displaystyle\sum_{l=1}^m 2 \Big(\int_{E_l} \frac{U^N_l - u^N(s)}{Z_l} ds - I_l \Big )\int_{E_l} \frac{1}{Z_l} (\delta_{lk} - w^N_k(s))ds\nonumber\\ + 2 \beta (U_k^N - U_k^*)\Big ], k=1,...,m\label{gradient_U}.
\end{gather}
\item[{ Step 7.}] Replace $(\tilde{A}^{N+1} , \tilde{U}^{N+1})$ with $(A^{N+1},U^{N+1}) \in V_R$ as follows
\begin{gather}
  a^{N+1}_{ij}(x)=\begin{cases}
\mu, &\quad{if} \  \tilde{a}^{N+1}_{ij}(x)\leq \mu,
\\
\tilde{a}^{N+1}_{ij}(x), &\quad{if} \ \mu\leq \tilde{a}^{N+1}_{ij}(x)\leq R,
\\
R, &\quad{if} \   \tilde{a}^{N+1}_{ij}(x)>R.
\end{cases}\label{gradientprojection_A}
  \\
  U^{N+1}_k=  \tilde{U}^{N+1}_k-\frac{1}{m}\sum_{k=1}^m \tilde{U}_k^{N+1}\label{gradientprojection_U}, \ k = 1,...,m
\end{gather}
Then replace $N$ with $N+1$ and move to Step 2.
\end{description}
Based on formula \eqref{eq:Frechet_Derivative1} similar algorithm is implemented for solving Problem $\mathcal{K}$.
\begin{remark}\label{EIT gradient}
Differentiability result and optimality condition similar to Theorem~\ref{thm:Frechet_Differentiability} and Corollary~\ref{optimalitycondition} are true for the Problem $\mathcal{I}$ and the gradient $\mathcal{I}_U'$ coincides
with $\mathcal{J}_U'$ from \eqref{eq:Frechet_Derivative}.
Similar algorithm for the gradient method in $\mathbb{R}^m$ applies to the Problem $\mathcal{I}$ in which case only iteration of the parameter $U$ is pursued.
\end{remark}

\section{Proofs of the Main Results }
\label{sec:Proofs_Main_Results}
Well-posedness of the elliptic problems \eqref{ general pde}--\eqref{bdry cond. of gen pde} and  \eqref{adj pde}--\eqref{bdry cond of adj pde}
follow from the Lax-Milgram theorem (\cite{evans1998partial}).
\begin{lemma}\label{existence lax}
For $\forall v\in V_R$ there exists a unique solution $u=u(\cdot,v)\in H^1(Q)$ to the problem  \eqref{ general pde}--\eqref{bdry cond. of gen pde}
which satisfy the energy estimate
\begin{align}
\norm{u}_{H^1(Q)}^2 \leq C \displaystyle\sum_{l=1}^{m} Z_l^{-2}U_l^2\label{Energy Estimate}
\end{align}
\end{lemma}
\textbf{Proof}: {\it Step 1. Introduction of the equivalent norm in $H^1(Q)$}. Let
\begin{equation}\label{newnorm}
|||u|||_{H^1(Q)}:= \bigg[ \int_{Q} |\nabla u|^2 dx +\displaystyle\sum_{l=1}^{m} \int_{E_l} u^2 ds \bigg ]^ \frac{1}{2},
\end{equation}
and prove that this is equivalent to the standard norm of $H^1(Q)$, i.e. there is $c>1$ such that
$\forall u\in H^1(Q)$
\begin{equation}
c^{-1}   \| u\|_{H^1(Q)}   \leq  |||u|||_{H^1(Q)}  \leq  c \|u\|_{H^1(Q)}
\label{equ-inequality}
\end{equation}
The second inequality immediately follows due to bounded embedding $H^1(Q)\hookrightarrow L^2(\partial Q)$ (\cite{evans1998partial}). To prove the first inequality assume on the contrary that
\begin{equation*}
\forall k>0 , \quad \exists u_k\in H^1(Q) \quad\text{such that} \ \|u_k\|_{H^1(Q)} >k |||u_k|||_{H^1(Q)}.
\end{equation*}
Without loss of generality we can assume that $\|u_k\| = 1$, and therefore
%\begin{equation*}
 %\||u_k|\|_{H^1(Q)} < \frac{1}{k}
%\end{equation*}
%This inequality implies:
\begin{equation}\label{limit}
\|\nabla u_k\|_{L_2(Q)} {\longrightarrow} \ 0, \ \|u_k\|_{L_2(E_l)} {\longrightarrow} \ 0, \quad\text{as} \ k \to \infty,  \quad l =1,2,...m.
\end{equation}
Since $\{ u_k \}$ is a bounded sequence in $H^1(Q)$, it is weakly precompact in $H^1(Q)$ and strongly precompact in both $L_2(Q)$ and $L_2(\partial Q)$ (\cite{nikolskii75,besov79,besov79a}). Therefore, there exists a subsequence $\{u_{k_j}\}$ and $u \in H^1(Q) $ such that $u_{k_j}$ converges to $u$ weakly in $H^1(Q)$ and strongly in $L_2(Q)$ and $L_2(\partial Q)$.
Without loss of generality we can assume that the whole sequence $\{ u_k \}$ converges to $u$. From the first relation of \eqref{limit} it follows that $\nabla u_k$ converges to zero strongly, and therefore also weakly in $L^2(Q)$. Due to uniqueness of the limit $\nabla u = 0$, and therefore $u=const$ a.e. in $Q$, and on the $\partial Q$ in the sense of traces. According to the second relation in \eqref{limit}, and since $|E_l| >0$, it follows that $const = 0$. This fact contradicts with $\|u_k\|=1$, and therefore the second inequality is proved.

{\it Step 2. Application of the Lax-Milgram theorem}.
Since $v\in V_R$, by using Cauchy-Bunyakowski-Schwartz (CBS) inequality, bounded trace embedding $H^1(Q)\hookrightarrow L^2(\partial Q)$ and \eqref{equ-inequality} we have the following estimations for the bilinear form $B$:
\begin{equation}
|B[u,\eta]|  \leq \alpha \|u\|_{H^1( Q)} \| \eta \|_{H^1( Q)}, \ B [u,u] \geq \beta \|u \|^2_{H^1(Q)}
\label{Lax 2 gen pde proof}
\end{equation}
where $\alpha, \beta >0$ are independent of $u,\eta$. Note that the component $U$ of the control vector $v$ defines a bounded linear functional $\hat{U}:H^1(Q)\to \mathbb{R}$ according to the right-hand side of \eqref{weaksolution}:
\begin{equation}
\hat{U}(\eta):= \displaystyle\sum_{l=1}^m \frac{U_l}{Z_l}\int_{E_l} \eta ds.
\label{righthandside}
\end{equation}
Indeed, by using CBS inequality and bounded trace embedding $H^1(Q)\hookrightarrow L^2(\partial Q)$ we have
\begin{equation}
| \hat{U}(\eta) | \leq |Q|^{\frac{1}{2}}    \bigl( \displaystyle\sum_{l=1}^m Z_l^{-2}U_l^2  \bigr)^{\frac{1}{2}} \|\eta\|_{L_2(\partial Q)}\leq C \|\eta\|_{H^1(Q)}
\label{F(eta) is bounded}
\end{equation}
From \eqref{Lax 2 gen pde proof},\eqref{F(eta) is bounded} and Lax-Milgram theorem (\cite{evans1998partial}) it follows that there exists a unique solution of the
problem \eqref{ general pde}--\eqref{bdry cond. of gen pde} in the sense of Definition~\ref{weaksolution}.

{\it Step 3. Energy estimate}.
By choosing $\eta$ as a weak solution $u$ in \eqref{weaksolution}, using \eqref{reconstructive. inequality.} and Cauchy's inequality with $\epsilon$ we derive
\begin{equation}
\mu \|\nabla u\|^2_{L_2(Q)}+ z_0 \displaystyle\sum_{l=1}^{m} \|u\|^2_{L_2(E_l)} \leq  \frac{c}{\epsilon} \displaystyle\sum_{l=1}^{m}  Z_l^{-2}U_l^2 + \epsilon |\partial Q| \displaystyle\sum_{l=1}^{m} \bigg( \int_{E_l} |u|^2 ds \bigg)
\label{eq:upper_bound_for_b[u,u]}
\end{equation}
where $z_0=\displaystyle\min_ {1\leq l \leq m} Z_l^{-1}$. By choosing $\epsilon = (2|\partial Q|)^{-1}z_0$ from \eqref{eq:upper_bound_for_b[u,u]} it follows that
\begin{equation}
|||u|||_{H^1(Q) } \leq C \displaystyle\sum_{l=1}^{m} Z_l^{-2}U_l^2.
\label{energy1}
\end{equation}
From \eqref{equ-inequality} and \eqref{energy1}, energy estimate \eqref{Energy Estimate} follows. Lemma is proved.$\blacksquare$
\begin{corollary}\label{adjoinedenergy}
For $\forall v\in V_R$ there exists a unique solution $\psi=\psi(\cdot,v)\in H^1(Q)$ of the adjoined problem \eqref{adj pde}--\eqref{bdry cond of adj pde}
which satisfy the energy estimate
\begin{align}
\| \psi \|_{H^1(Q)}^2 \leq C \displaystyle\sum_{l=1}^{m} Z_l^{-2}\Big[ \int_{E_l} \frac{U_l-u}{Z_l} ds -  I_l\Big]^2
\label{Adj Energy Estimate}
\end{align}
where $u=u(\cdot;v)\in H^1(Q)$ is a solution of the problem \eqref{ general pde}--\eqref{bdry cond. of gen pde} .
\end{corollary}

\textbf{Proof of Theorem~\ref{thm:continuity_of_J}}.
Let $\{ v_k\} = \{ (A^k, U^k)\}\subset V_R$ be a minimizing sequence
\begin{equation*}
\lim_{k \to \infty} \mathcal{J}(v_k) = \mathcal{J}_*
\end{equation*}
Since $\{A^k\}$ is a bounded sequence in $H^{\epsilon}(Q;\mathbb{M}^{n \times n})$, it is weakly precompact in $H^{\epsilon}(Q;\mathbb{M}^{n \times n})$ and strongly  precompact in $L_2(Q;\mathbb{M}^{n \times n})$ (\cite{nikolskii75,besov79,besov79a}). Therefore, there exists a subsequence  $\{ A^{k_p} \}$ which converges weakly in $H^{\epsilon}(Q;\mathbb{M}^{n \times n})$ and strongly in $L_2(Q;\mathbb{M}^{n \times n})$ to some element $A \in H^{\epsilon}(Q;\mathbb{M}^{n \times n})$. Since any strong convergent sequence in $L_2(Q;\mathbb{M}^{n \times n})$ has a subsequence which converges a.e. in $Q$, without loss of generality one can assume that the subsequence $A^{k_p}$ converges to $A$ a.e. in $Q$, which implies that $A\in L_\infty(Q;\mathbb{M}^{n \times n})\cap H^{\epsilon}(Q;\mathbb{M}^{n \times n})\cap V_R$. Since $U^k$ is a bounded sequence in $\mathbb{R}^m$ it has a subsequence which converges to some $U\in\mathbb{R}^m, |U|\leq R$. Without loss of generality we cam assume that the whole minimizing sequence $v_k = (A_k, U^k)$ converges $v = (A, U)\in V_R$ in the indicated way.

Let $u_k= u(x;v_k)$, $u= u(x;v)\in H^1(Q)$ are weak solutions of \eqref{ general pde}--\eqref{bdry cond. of gen pde} corresponding to $v_k$ and $v$ respectively. By Lemma~ \ref{existence lax} $u_k$ satisfy the energy estimate \eqref{Energy Estimate} with $U^k$ on the right hand side, and therefore it is uniformly bounded in $H^1(Q)$. By the Rellich-Kondrachov compact embedding theorem there exists a subsequence  $\{u_{k_p} \}$ which converges weakly in $H^1(Q)$ and strongly in both $L_2(Q)$ and $L_2(\partial Q)$ to some function $\tilde{u} \in H^1(Q)$(\cite{nikolskii75,besov79,besov79a}). Without loss of generality assume that the whole sequence $u_k$ converges to $\tilde{u}$ weakly in $H^1(Q)$ and strongly both in $L_2(Q)$ and $L_2(\partial Q)$. For any fixed $\eta \in C^1(Q)$ weak solution $u_k$ satisfies the following integral identity
\begin{equation}\label{identityforuk}
\int_Q  \displaystyle\sum_{i,j=1}^n a^k_{ij}u_{k_{x_j}} \eta_{x_i}  dx + \displaystyle\sum_{l=1}^m \frac{1}{Z_l}\int_{E_l} u_k \eta ds  = \displaystyle\sum_{l=1}^m \frac{1}{Z_l}\int_{E_l} \eta  U^k_l ds.
\end{equation}
Due to weak convergence of $\nabla u_k$ to $\nabla \tilde{u}$ in $L_2(Q;\mathbb{R}^n)$, strong convergence of $u_k$ to $\tilde{u}$ in $L_2(\partial Q)$, strong convergence of $a^k_{ij}$ to $a_{ij}$ in $L_2(Q)$ and convergence of $U^k$ to $U$, passing to the limit as $k\to \infty$, from \eqref{identityforuk} it follows
\begin{equation}\label{identityforu}
\int_Q  \displaystyle\sum_{i,j=1}^n a_{ij}\tilde{u}_{{x_j}} \eta_{x_i}  dx + \displaystyle\sum_{l=1}^m \frac{1}{Z_l}\int_{E_l} \tilde{u} \eta ds  = \displaystyle\sum_{l=1}^m \frac{1}{Z_l}\int_{E_l} \eta  U_l ds.
\end{equation}
Due to density of $C^1(Q)$ in $H^1(Q)$ (\cite{nikolskii75,besov79,besov79a}) the integral identity \eqref{identityforu} is true for arbitrary $\eta \in H^1(Q)$. Hence, $\tilde{u}$ is a weak solution of the problem \eqref{ general pde}--\eqref{bdry cond. of gen pde} corresponding to the control vector $v=(A,U)\in V_R$. Due to uniqueness of the weak solution it follows that
$\tilde{u}=u$, and the sequence $u_k$ converges to the weak solution $u=u(x;v)$ weakly in $H^1(Q)$, and strongly both in $L_2(Q)$ and $L_2(\partial Q)$. The latter easily implies that
\begin{equation*}
\mathcal{J}(v) = \lim_{n \to \infty} \mathcal{J}(v_n) = \mathcal{J}_*
\end{equation*}
Therefore, $v\in V_*$ is an optimal control and \eqref{solutionexists} is proved. $\blacksquare$

\textbf{Proof of Theorem}~\ref{thm:Frechet_Differentiability}.
Let $v=(A,U)\in V_R$ is fixed and $\delta v=(\delta A, \delta U)$ is an increment such that $\bar{v}=v+\delta v \in V_R$ and
$u=u(\cdot;v), \bar{u}=u(\cdot;v+\delta v) \in H^1(Q)$ are respective weak solutions of the problem \eqref{ general pde}--\eqref{bdry cond. of gen pde}.
Since $u(\cdot;A,U)$ is a linear function of $U$ it easily follows that
\[ w^k = \frac{\partial u}{\partial U_k}=u(\cdot;A,e_k)\in H^1(Q), \ k = 1,2,..,m \]
 is a solution of \eqref{ general pde}--\eqref{bdry cond. of gen pde} with $v=(A,e_k) $, $e_k\in\mathbb{R}^m$ is a unit ort vector in $x_k$-direction.
 Straightforward calculation imply that
 \begin{align*}
\frac{\partial \mathcal{J}}{\partial U_k} = \displaystyle\sum_{l=1}^m 2 \Big [\int_{E_l} \frac{U_l - u}{Z_l} ds - I_l\Big ]  \int_{E_l} \frac{1}{Z_l} ( \delta_{lk} - w^k ) ds  + 2 \beta (U_k - U_k^*), \ k=1,...,m.
\end{align*}
where $\delta_{lk}$ is a Kronecker delta.

In order to prove the Fr\'echet differentiability with respect to $A$, assume that $\delta U=0$ and transform the increment of  $\mathcal{J}$ as follows
\begin{gather}
\delta \mathcal{J}:= \mathcal{J}(v+\delta v)-\mathcal{J}(v)  =  \displaystyle\sum _{l=1}^m  \frac{1}{Z_l}\int_{E_l}  2 \bigg(\int_{E_l} \frac{u-U_l}{Z_l} ds +  I_l\bigg) \delta u ds + R_1,
\label{variation of obj}\\
R_1= \displaystyle\sum _{l=1}^m Z_l^{-2}\bigg(\int_{E_l} \delta u ds\bigg)^2 \leq \displaystyle\sum _{l=1}^m |E_l|Z_l^{-2} |||\delta u |||^2_{H^1(Q)},
\label{quadraticterm}
\end{gather}
where $\delta u = \bar{u}-u$. By subtracting integral identities \eqref{weaksolution} for $\bar{u}$ and $u$, and by choosing test function $\eta=\psi(\cdot;v)$ as a solution of the adjoined problem \eqref{adj pde}--\eqref{bdry cond of adj pde} we have
\begin{equation}\label{gradient1}
\displaystyle\int_{Q}  \displaystyle\sum_{ij} \bigg( \delta a_{ij} u_{x_j} + a_{ij} (\delta u)_{x_j}  +  \delta a_{ij} (\delta u)_{x_j}     \bigg) \psi_{x_i}dx
+  \displaystyle\sum_{l=1}^m \frac{1}{Z_l}\int_{E_l}  \psi \delta u  ds = 0.
\end{equation}
By choosing $\eta=\delta u$ in the integral identity \eqref{weaksolutionadjoined} for the weak solution $\psi$ of the adjoined problem we have
\begin{equation}\label{gradient2}
-\displaystyle\int_Q \displaystyle\sum_{ij} a_{ij} \psi_{x_i} \delta u_{x_j} dx  + \displaystyle\sum_{l} \int_{E_l} \frac{\delta u}{Z_l} \big[  2\int_{E_l} \frac{u-U_l}{Z_l} ds + 2 I_l - \psi  \big] ds=0
\end{equation}
Adding \eqref{gradient1} and \eqref{gradient2} we derive
\begin{gather}
\displaystyle\sum _{l=1}^m \displaystyle \frac{1}{Z_l} \int_{E_l} 2\bigg(  \displaystyle\int_{E_l} \frac{u-U_l}{Z_l}dS(x) +I_l \bigg) \delta u ds
=  \displaystyle\int_Q  (- \displaystyle\sum_{ij}   \delta a_{ij} u_{x_j}  \psi_{x_i} -\displaystyle\sum_{ij}   \delta a_{ij} (\delta u)_{x_j}  \psi_{x_i}) dx.
\label{gradient3}
\end{gather}
From \eqref{variation of obj} and \eqref{gradient3} it follows that
\begin{equation}\label{gradient4}
\delta \mathcal{J} = -\displaystyle\int_Q   \displaystyle\sum_{ij}  u_{x_j}  \psi_{x_i}  \delta a_{ij} dx + R_1+R_2
\end{equation}
where
\begin{equation}\label{remainder}
R_2= - \displaystyle\int_Q  \displaystyle\sum_{ij}   \delta a_{ij} (\delta u)_{x_j}  \psi_{x_i} dx.
\end{equation}
To complete the proof it remains to prove that
\begin{equation}\label{osmall}
R_1+R_2 = o(\| \delta A \|_{L_\infty(Q;M^{n\times n})})  \quad\text{as} \  \| \delta A \|_{L_\infty(Q;M^{n\times n}} \to 0.
\end{equation}
By subtracting integral identities \eqref{weaksolution} for $\bar{u}$ and $u$ again, and by choosing test function $\eta=\delta u$ we have
\begin{equation}\label{osmall1}
\displaystyle\int_Q \displaystyle\sum_{ij} \bar{a}_{ij} (\delta u)_{x_j} (\delta u)_{x_i} dx +  \displaystyle\sum_{l=1}^m \frac{1}{Z_l}\int_{E_l}  (\delta u)^2  ds =
-\displaystyle\int_{Q}  \displaystyle\sum_{ij} \delta a_{ij} u_{x_j} (\delta u)_{x_i} dx.
\end{equation}
By using positive definiteness of  $\bar{A}\in V_R$ and by applying Cauchy inequality with $\epsilon >0$ to the right hand side, from \eqref{osmall1} it follows that
\begin{equation}
 \mu \displaystyle\int_{Q}  |\nabla \delta u|^2 dx +  \displaystyle\sum_{l=1}^m \frac{1}{Z_l}\int_{E_l}  (\delta u)^2 ds   \leq \epsilon \displaystyle\int_{Q}  |\nabla \delta u|^2 dx + \frac{c}{\epsilon} \int_{Q}   |\displaystyle\sum_{ij} \delta a_{ij}|^2  |\nabla u|^2.
\label{Energy est delta u step 1}
\end{equation}
By choosing $\epsilon=\mu/2$ and by applying the energy estimate \eqref{Energy Estimate} from \eqref{Energy est delta u step 1} we derive
\begin{equation}
|||\delta u|||^2_{H^1(Q)} \leq  C  \| \delta A\|^2_{L_{\infty} (Q;\mathbb{M}^{n\times n})}.
\label{energy estimate for delta u}
\end{equation}
From \eqref{remainder} it follows that
\begin{equation}\label{remainder1}
|R_2|\leq C  \| \delta A\|_{L_{\infty} (Q;\mathbb{M}^{n\times n})} \| \nabla \delta u \|_{L_2(Q)} \| \nabla \psi \|_{L_2(Q)}.
\end{equation}
From \eqref{Energy Estimate}, \eqref{equ-inequality}, \eqref{Adj Energy Estimate}, \eqref{quadraticterm}, \eqref{energy estimate for delta u} and \eqref{remainder1}, desired estimation \eqref{osmall} follows. Theorem is proved.$\blacksquare$
\begin{remark}\label{convex}Functional \eqref{eq:cost_functional1} in the optimal control Problem $\mI$ is convex due to the following formula
\begin{equation*}
{\mathcal I}(\alpha U^1 +(1-\alpha)U^2)=\alpha {\mI}(U^1)+(1-\alpha) {\mI}(U^2)-\alpha(1-\alpha)\displaystyle\sum_{l=1}^m Z_l^{-2} \Big |\displaystyle\int_{E_l}(U^1_l-U^2_l -u^1+u^2)ds\Big |^2
\end{equation*}
where $U^1,U^2\in W, \alpha\in[0,1]; u^i=u(\cdot;U^i), i=1,2$ is a solution of \eqref{ general pde}--\eqref{bdry cond. of gen pde} with $U=U^i$. Therefore, unique solution of the EIT problem would be a unique global minimizer of the Problem ${\mathcal I}$.
\end{remark}
\section{Numerical Results}
\label{sec:numerical_results}

In this section we describe computational results for solving the Inverse EIT Problem
in the 2D case ($n=2$) according to the algorithm outlined
in Section~\ref{sec:Gradient Method in Banach Space}. First, we discuss the structure of our 2D computational
model. The complexity level of this model is chosen to adequately represent the diagnosis of breast cancer in reality.
Then we briefly describe the numerical approaches
used for discretizing the problem in space and accurately solving related PDEs to ensure advanced
performance of numerical techniques included in the optimization framework, e.g.~PCA-based
re-parameterization and proper regularization. Finally, we show the outcomes
of applying the proposed computational algorithm to this 2D model and discuss further steps to improve
the performance.

\subsection{Computational Model in 2D Space}
\label{sec:model}
\begin{figure}[!htb]
  \begin{center}
  \mbox{
  \subfigure[]{\includegraphics[width=0.33\textwidth]{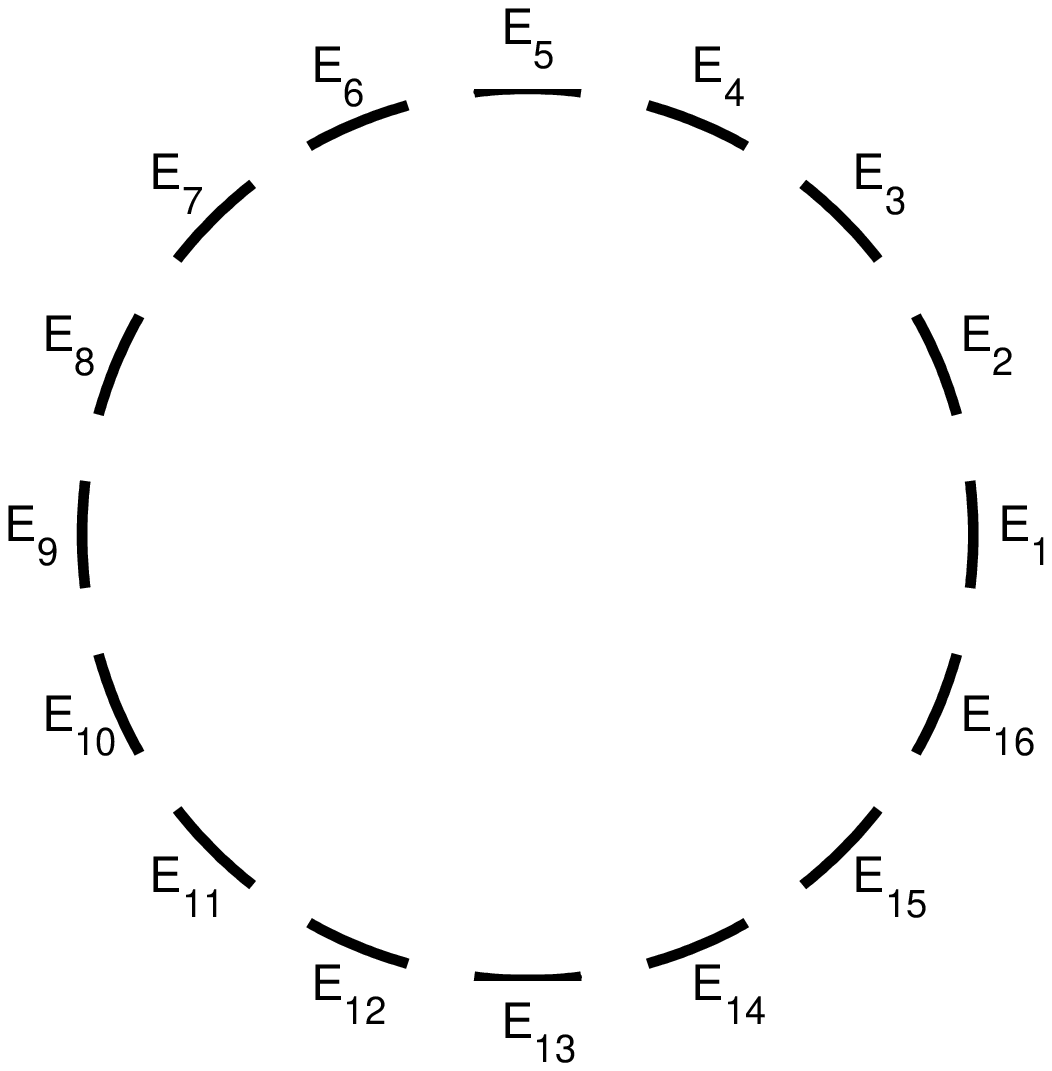}}
  \subfigure[]{\includegraphics[width=0.33\textwidth]{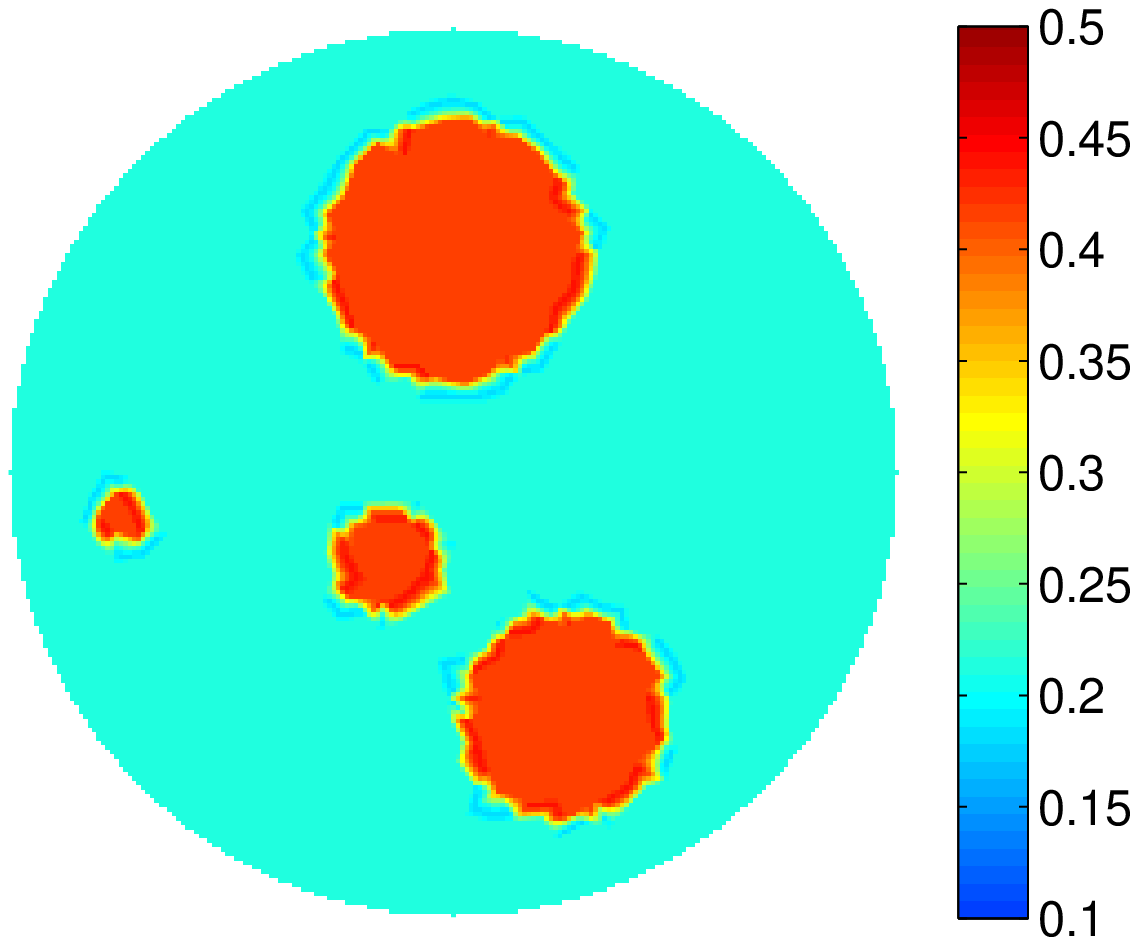}}
  \subfigure[]{\includegraphics[width=0.33\textwidth]{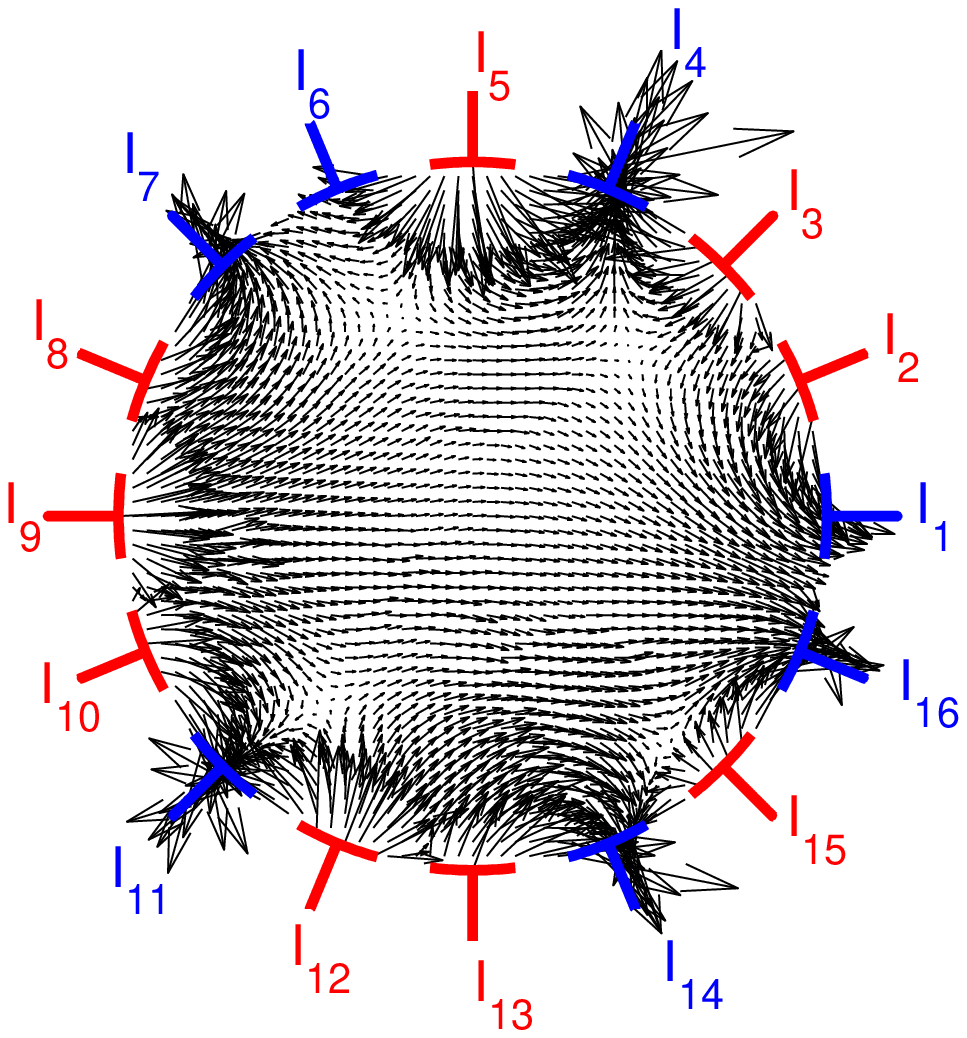}}}
  \caption{(a) Equispaced geometry of electrodes $E_l$ placed over boundary $\partial Q$.
    (b) True electrical conductivity $\sigma_{true}(x)$.
    (c) Electrical currents $I_l$ (positive in red, negative in blue) injected by electrodes $E_l$.
    Black arrows show the distribution of flux $\sigma(x) \nabla u(x)$ of the electrical potential $u$
    in the interior of domain $Q$.}
  \label{fig:geometry}
  \end{center}
\end{figure}

We pursue computational analysis of the inverse EIT problem with removed assumption on anisotropy for electrical
conductivity tensor $A(x)$, i.e.~$A(x) = \sigma(x) I$, where $I$ is a $2\times 2$ unit matrix.
Problem ${\mathcal J}$ consists of the minimization of cost functional $\mathcal{J}(\sigma, U)$
defined in \eqref{eq:cost_functional} on control set $V_R$, where $u=u(\cdot; \sigma, U)$ solves the elliptic PDE problem
\begin{alignat}{3}
  & {\rm div} \left( \sigma(x) \nabla u(x) \right) = 0, & \quad  x \in Q \label{eq:forward_1}\\
  & \Dpartial{u(x)}{n} = 0, & \quad x \in \partial Q - \bigcup\limits_{l=1}^{m} E_l \label{eq:forward_2}\\
  & u(x) +  Z_l \sigma(x) \Dpartial{u(x)}{n} = U_{l}, & \quad x \in E_l, \ l= \overline{1,m} \label{eq:forward_3}
  %& \displaystyle\int_{E_l} \sigma(x) \Dpartial{u(x)}{n} ds = I_l, & \quad x\in E_l, \ l= \overline{1,m} \label{eq:forward_4}
\end{alignat}
where $n$ is an external unit normal vector on $\partial Q$. The first term in the cost functional $\mathcal{J}(\sigma,U)$ characterizes
mismatch of the condition
\begin{equation}\label{eq:forward_4}
\displaystyle\int_{E_l}  \sigma(x)\Dpartial{u(x)}{n} ds = I_l, \quad l= \overline{1,m}
\end{equation}
in light of the Robin condition \eqref{eq:forward_3}
We choose $Q$ as a disk
\begin{equation}
  Q = \left\{ x \in \RR^{2}: \ x_1^2 + x_2^2 < r_Q^2 \right\}
  \label{eq:domain2D}
\end{equation}
of radius $r_Q = 0.1$ with $m = 16$ equidistant electrodes $E_l$ with half-width $w = 0.12$ rad covering approximately
61\% of boundary $\partial Q$ as shown in Figure \ref{fig:geometry}(a).

The actual (true) electrical conductivity $\sigma_{true}(x)$ we seek to reconstruct is given analytically by
\begin{equation}
  \sigma_{true}(x)= \left\{
  \begin{aligned}
  & 0.4, \qquad &x_1^2+(x_2-0.05)^2 \leq (0.03)^2\\
  & 0.4, \qquad &(x_1+0.075)^2+(x_2+0.01)^2 \leq (0.0063)^2\\
  & 0.4, \qquad &(x_1+0.015)^2+(x_2+0.02)^2 \leq (0.0122)^2\\
  & 0.4, \qquad &(x_1-0.025)^2+(x_2+0.055)^2 \leq (0.0235)^2\\
  & 0.2, \qquad &\text{otherwise}
  \end{aligned} \right.
  \label{eq:sigma_true}
\end{equation}
measured in $(Ohm \cdot m)^{-1}$ and setting $\sigma_c = 0.4$ for cancer-affected parts (4 spots of different size)
and $\sigma_h = 0.2$ to healthy tissues parts as seen in Figure~\ref{fig:geometry}(b). Electrical currents $I_l$
injected by electrodes $E_l$ are provided in Table~\ref{tab:curr_model_C} and shown schematically in Figure~\ref{fig:geometry}(c).
This figure also shows the distribution of flux $\sigma(x) \nabla u(x)$
of the electrical potential $u$ in the interior of domain $Q$ corresponding to $\sigma_{true}(x)$.
\begin{table}[!htb]
  \caption{``Current--to--voltage'' model parameters: electrical currents $I_l$ injected by electrodes $E_l, \, l = 1, \ldots, 16$,
    with contact impedances $Z_l$, and initial guess for boundary voltages $U_{l,ini}$. The unit system used for all values
    is SI.}
  \label{tab:curr_model_C}
  \begin{center}
  \begin{tabular}{lcccccccccccccccc}
  \hline\noalign{\smallskip}
  Electrode, $l$ & 1 & 2 & 3 & 4 & 5 & 6 & 7 & 8 & 9 & 10 & 11 & 12 & 13 & 14 & 15 & 16\\
  \noalign{\smallskip}\hline\noalign{\smallskip}
  $I_l \cdot 10^2$, A & -3 & 2 & 3 & -7 & 6 & -1 & -4 & 2 & 4 & 3 & -5 & 4 & 3 & -5 & 2 & -4\\
  $Z_l \cdot 10^1$, Ohm & 1 & 1 & 1 & 1 & 1 & 1 & 1 & 1 & 1 & 1 & 1 & 1 & 1 & 1 & 1 & 1\\
  \noalign{\smallskip}\hline\noalign{\smallskip}
  $U_{l,ini}$, V & -1 & 1 & -1 & 1 & -1 & 1 & -1 & 1 & -1 & 1 & -1 & 1 & -1 & 1 & -1 & 1\\
  \noalign{\smallskip}\hline
  \end{tabular}
  \end{center}
\end{table}

Our optimization framework integrates computational facilities for solving state PDE problem
\eqref{eq:forward_1}--\eqref{eq:forward_3}, adjoint PDE problem \eqref{adj pde}--\eqref{bdry cond of adj pde},
and evaluation of the Fr\'echet gradient according to \eqref{eq:Frechet_Derivative}, \eqref{eq:Frechet_Derivative1}. These facilities are incorporated
by using {\tt FreeFem++}, see \cite{FreeFem} for details, an open--source, high--level integrated
development environment for obtaining numerical solutions of PDEs based on the Finite Element Method.
To numerically solve the state PDE problem \eqref{eq:forward_1}--\eqref{eq:forward_3}, spatial discretization
is carried out by implementing triangular finite elements, P2 piecewise quadratic (continuous)
representation for electrical potential $u(x)$ and P0 piecewise constant representation for conductivity field
$\sigma(x)$. The system of algebraic equations obtained after such discretization is solved with
{\tt UMFPACK}, a solver for nonsymmetric sparse linear systems. The same technique is used for the numerical
solution of adjoint problem \eqref{adj pde}--\eqref{bdry cond of adj pde}. All computations are performed
using 2D domain $Q$ \eqref{eq:domain2D} which is discretized using mesh $\mM(n_v)$ created by specifying
$n_v = 96$ vertices over boundary $\partial Q$ and totaling 1996 triangular finite elements inside $Q$.

In terms of the initial guess in the iterative algorithm shown in Section~\ref{sec:Main_Results}, unless stated otherwise,
we take a constant approximation to \eqref{eq:sigma_true}, given by
$\sigma_{ini} = \frac{1}{2} \left( \sigma_h + \sigma_c \right) = 0.3$. Initial guess for boundary voltages is provided
in Table~\ref{tab:curr_model_C} which is consistent with the ground potential condition \eqref{grounding}.
Determining the Robin part of the boundary conditions in \eqref{eq:forward_3} we equally set the electrode
contact impedance $Z_l = 0.1$.

The iterative optimization algorithm is performed by the Sparse Nonlinear OPTimizer~{\tt SNOPT}, a software package
for solving large-scale nonlinear optimization problems, see \cite{GillMurraySaunders02}. It employs a sparse sequential
quadratic programming~(SQP) algorithm with limited-memory quasi-Newton approximations to the Hessian of the Lagrangian.
This makes {\tt SNOPT} especially effective for nonlinear problems with computationally expensive functionals and gradients,
like in our problem. The termination conditions set for {\tt SNOPT} are
$\left| \frac{\J^{N}-\J^{N-1}}{\J^{N-1}} \right| < 10^{-6}$ or maximum number of optimization iterations
$N_{max} = 250$ whichever comes first.

\subsection{Reduced--Dimensional Optimization via PCA--based Re-parameterization}
\label{sec:PCA}

From a viewpoint of numerical optimization, Problem ${\mathcal J}$ in its spatially discretized form is
over-parameterized even for moderate size models. As previously mentioned in Section~\ref{sec:model}, our
2D computational model requires a solution for 1996-component electrical conductivity vector $\sigma$ when using
relatively coarse mesh $\mM(96)$. %In fact, the amount of observed data ($m = 16$) is very small compared to
To overcome ill-posedness due to over-parameterization we implement
re-parameterization of the control set based on PCA, which is also known as Proper Orthogonal Decomposition (POD)
or Karhunen--Lo\`{e}ve Expansion.

Without loss of generality, we consider a model which contains $N_{\sigma}$ model parameters. We assume the existence
of a set of $N_r$ sample solutions (realizations) $\sigma_j$, $j = 1, \dots, N_r$, each of size $N_{\sigma}$.
%(see examples in Figure~\ref{fig:PCA samples}).
For simplicity we assume a Gaussian (normal) distribution for the model parameters,
i.e., $\sigma \sim \NN(\bar \sigma, \CC_M)$, where $\bar \sigma = (1/N_r) \sum_{j=1}^{N_r} \sigma_j$. Covariance matrix
$\CC_M$ may be approximated by
\begin{equation}
  \CC_M \approx \dfrac{X X^T}{N_r-1}, \quad
  X^{N_{\sigma} \times N_r} =
  \left[ \sigma_1 - \bar \sigma \ \ \dots \ \  \sigma_{N_r} - \bar \sigma\right].
  \label{eq:cov_mat}
\end{equation}
It is more efficient to perform singular value decomposition (SVD) on matrix $Y = X / \sqrt{N_r-1}$ of size
$N_{\sigma} \times N_r$, rather than on covariance matrix $\CC_M$ of size $N_{\sigma} \times N_{\sigma}$, as
$N_r \ll N_{\sigma}$. The SVD factorization with truncation is then applied to matrix $Y$
\begin{equation}
  Y \approx \tilde U_{N_{\xi}} \, \tilde \Sigma_{N_{\xi}} \, \tilde V^T_{N_{\xi}},
\end{equation}
where diagonal matrix $\tilde \Sigma_{N_{\xi}}$ contains the singular values of $Y$, and matrices $\tilde U_{N_{\xi}}$
and $\tilde V^T_{N_{\xi}}$ are matrices containing the left and right singular vectors of $Y$. More specifically,
matrix $\tilde \Sigma_{N_{\xi}}$ is truncated to keep only $N_{\xi}$ singular values. Similarly, analogous truncations are
applied to $\tilde U_{N_{\xi}}$ and $\tilde V^T_{N_{\xi}}$.

We define a linear transformation
\begin{equation}
  \Phi^{N_{\sigma} \times N_{\xi}} =
  \tilde U_{N_{\xi}} \, \tilde \Sigma_{N_{\xi}}, \qquad
  N_{\xi} \leq N_{\rm min} = \min \{ N_{\sigma}, N_r \}
  \label{eq:Phi_trunc}
\end{equation}
to project the initial control space defined for model parameters $\sigma$ onto reduced-dimensional $\xi$-space
which contains only $N_{\xi}$ largest principal components, see \cite{bukshtynov15}, by means of the unique mapping
\begin{equation}
  \sigma = \Phi \, \xi + \bar \sigma.
  \label{eq:PCA_ctrls}
\end{equation}
To construct a ``backward'' mapping, the simplest approach is to approximate the inverse of matrix $\Phi$,
which cannot be inverted due to its size $N_{\sigma} \times N_{\xi}$, using a pseudo-inverse matrix $\hat \Phi^{-1}$
\begin{equation}
  \xi = \Phi^{-1} \, (\sigma - \bar \sigma) \approx \hat \Phi^{-1} \, (\sigma - \bar \sigma)
  = \tilde \Sigma^{-1}_{N_{\xi}} \, \tilde U^T_{N_{\xi}} \, (\sigma - \bar \sigma).
  \label{eq:PCA_backward}
\end{equation}

The optimal control problem defined in Section~\ref{sec:Optimal_Control_Problem} can now be restated in terms of
new model parameters $\xi$ used in place of control $\sigma$ as follows
\begin{equation}
  \left( \hat \xi, \hat U \right) = \underset{\xi,U}{\argmin}\,\mJ(\xi,U),
  \label{eq:minJ_xi}
\end{equation}
subject to discretized PDE model \eqref{eq:forward_1}--\eqref{eq:forward_3}, and using mappings given by
\eqref{eq:PCA_ctrls}--\eqref{eq:PCA_backward}. By applying \eqref{eq:PCA_ctrls} and the chain rule for derivatives,
gradient $\nabla_{\xi} \mJ$ of cost functional $\mJ$ with respect to controls $\xi$ can be expressed as
\begin{equation}
  \nabla_{\xi} \mJ = \Phi^T \, \nabla_{\sigma} \mJ =
  \tilde \Sigma_{N_{\xi}} \, \tilde U^T_{N_{\xi}} \, \nabla_{\sigma} \mJ.
  \label{eq:PCA_grad}
\end{equation}
This expression, in fact, defines projection of gradient $\nabla_{\sigma} \mJ$ shown in \eqref{eq:Frechet_Derivative}
from initial (physical) $\sigma$-space onto the reduced-dimensional $\xi$-space. A summary of the discretized finite-dimensional version
of the projective gradient method in Besov spaces for the Problem $\mathcal{J}$ outlined in Section~\ref{sec:Gradient Method in Banach Space}, employing solution of the problem~\eqref{eq:minJ_xi}, is provided in Algorithm~\ref{alg:param_opt}. The same algorithm could be easily adjusted for solving
Problem $\mathcal{I}$ in which case only iteration of $U$ is pursued (see Remarks~\ref{EIT gradient} and~\ref{convex}).
%
%\begin{figure}[!htb]
%  \begin{center}
%  \mbox{
%  %\subfigure{\includegraphics[width=0.33\textwidth]{Figs/PCA_1}}
%  \subfigure{\includegraphics[width=0.33\textwidth]{Figs/PCA_2}}
%  \subfigure{\includegraphics[width=0.33\textwidth]{Figs/PCA_3}}
%  \subfigure{\includegraphics[width=0.33\textwidth]{Figs/PCA_4}}}
%  \caption{Three random realizations from the set of $N_r = 500$ sample solutions $\sigma_j$, $j = 1, \dots, N_r$, created
%    to perform PCA-based re-parameterization for control vector $\sigma$.}
%  \label{fig:PCA samples}
%  \end{center}
%\end{figure}
%

A problem of approximating covariance matrix $\CC_M$ in \eqref{eq:cov_mat} to support our current 2D computational
model described in Section~\ref{sec:model} is solved in the following way. A set of $N_r = 500$ realizations $\sigma_j$,
$j = 1, \dots, N_r$, is created using a generator of (uniformly distributed) random numbers. Each realization ``contains''
from 1 to 7 ``cancer-affected'' areas with $\sigma_c = 0.4$. Each area is located randomly within domain $Q$ and represented
by a circle of a randomly chosen radius $0 < r \leq 0.3 r_Q$. %Three randomly chosen realizations from this set are shown in
%Figure~\ref{fig:PCA samples}.
We refer the discussion on choosing optimal number of principal components $N_{\xi}$ to Appendix~\ref{sec:tuning_pca}
to consider it as a part of a tuning process in optimizing the overall performance of our computational framework.

\begin{algorithm}[h]
  \begin{algorithmic}
    \STATE $N \leftarrow 0$
    \STATE $U^0 \leftarrow $ initial guess $U_{ini}$
    \STATE $\sigma^0 \leftarrow $ initial guess $\sigma_{ini}$
    \STATE construct $\Phi$ and $\hat \Phi^{-1}$ by \eqref{eq:Phi_trunc} and \eqref{eq:PCA_backward}
    \STATE $\xi^0 \leftarrow \sigma^0$ using \eqref{eq:PCA_backward}
    \REPEAT
    \STATE given estimate of $(\sigma^N, U^N)$, solve state equations \eqref{eq:forward_1}--\eqref{eq:forward_3} for $u^N$
    \STATE given $u^N$ and $(\sigma^N, U^N)$, solve adjoint equations \eqref{adj pde}--\eqref{bdry cond of adj pde} for $\psi^N$
    \STATE given estimate of $\sigma^N$, solve \eqref{eq:forward_1}--\eqref{eq:forward_3} for $v_k^N$ where $U = e_k$
    \STATE $\left( \nabla_{\sigma} \mJ(\sigma^N, U^N), \nabla_U \mJ(\sigma^N, U^N) \right)
      \leftarrow \sigma^N, U^N, u^N, \psi^N$ by \eqref{eq:Frechet_Derivative}
    \STATE $\nabla_{\xi} \mJ(\xi^N, U^N) \leftarrow \nabla_{\sigma} \mJ(\sigma^N, U^N)$ by \eqref{eq:PCA_grad}
    \STATE update $\xi^{N+1}$ and $U^{N+1}$ by using descent directions $D_{\xi}$ and $D_U$ obtained from
      $\nabla_{\xi} \mJ^N$ and $\nabla_U \mJ^N$:
      \begin{equation}
        \xi^{N+1} = \xi^N - \tau^N D_{\xi} \left( \nabla_{\xi} \mJ(\xi^N, U^N) \right)
      \end{equation}
      \begin{equation}
        U^{N+1} = U^N - \tau^N D_U \left( \nabla_U \mJ(\xi^N, U^N) \right)
      \end{equation}
    \STATE $\sigma^{N+1} \leftarrow \xi^{N+1}$ by \eqref{eq:PCA_ctrls}
    \STATE $N \leftarrow N + 1$
    \UNTIL termination criteria are satisfied to a given tolerance
  \end{algorithmic}
  \caption{Optimization workflow utilizing PCA-based control space re-parameterization}
  \label{alg:param_opt}
\end{algorithm}
\begin{remark} Corollary~\ref{cor:Frechet_Diffrentiability_K} in the context of the model example claims that the Fr\'echet gradient $\mathcal{K}'( \sigma , U)\in \textbf{ba} (Q) \times \mathbb{R}^m$ is
\begin{gather}
\mathcal{K}'(\sigma, U ) = \Big(\mathcal{K}'_\sigma (\sigma, U ), \mathcal{K}'_U(\sigma, U) \Big)=
\nonumber\\
\left( -\displaystyle\sum_{j=1}^m \nabla u_j \cdot \nabla\psi_j,  \Big( \displaystyle\sum_{j=1}^m \displaystyle\sum_{l=1}^m 2 \Big[\int_{E_l} \frac{U^j_l - u_j}{Z_l} ds - I^j_l\Big]  \int_{E_l} \frac{\delta_{l,\theta_{kj}} - w^{\theta_{kj}}(s)}{Z_l}ds + 2 \beta (U_k - U_k^*)  \Big)_{k=1}^m  \right) \label{eq:Frechet_Derivative2}
\end{gather}
\end{remark}

\subsection{Numerical Results for EIT and Inverse EIT Problems}
\label{sec:EIT_forward_backward}
To test the effectiveness of our gradient descent method, we simulate a realistic model example of the inverse EIT problem which adequately represent the diagnosis of the breast cancer
in reality. Simulation and computational analysis consists of three stages.
%The goal of this section is to pursue the development twofold:

{\it Stage 1.} By selecting boundary current pattern $I=(I_l)_{l=1}^{16}$ we simulate EIT model example with $\sigma=\sigma_{true}$ %U_{true} \right)$
by solving Problem $\mathcal{I}$ by the gradient descent
method described in Section~\ref{sec:Gradient Method in Banach Space}, Algorithm~\ref{alg:param_opt} and identifying optimal control $U_{true}$. 
Practical analogy of this step is implementation of the ``current--to--voltage''
procedure: by injecting current pattern $I=(I_l)_{l=1}^{16}$ on the electrodes $E_l, \ l= 1, \cdots, 16$, take the measurement
of the voltages $U^*=(U^*_1,\cdots,U^*_{16})$. In our numerical simulations $U_{true}$ is identified with $U^*$.

\begin{figure}[!htb]
  \begin{center}
  \mbox{
  \subfigure[]{\includegraphics[width=0.50\textwidth]{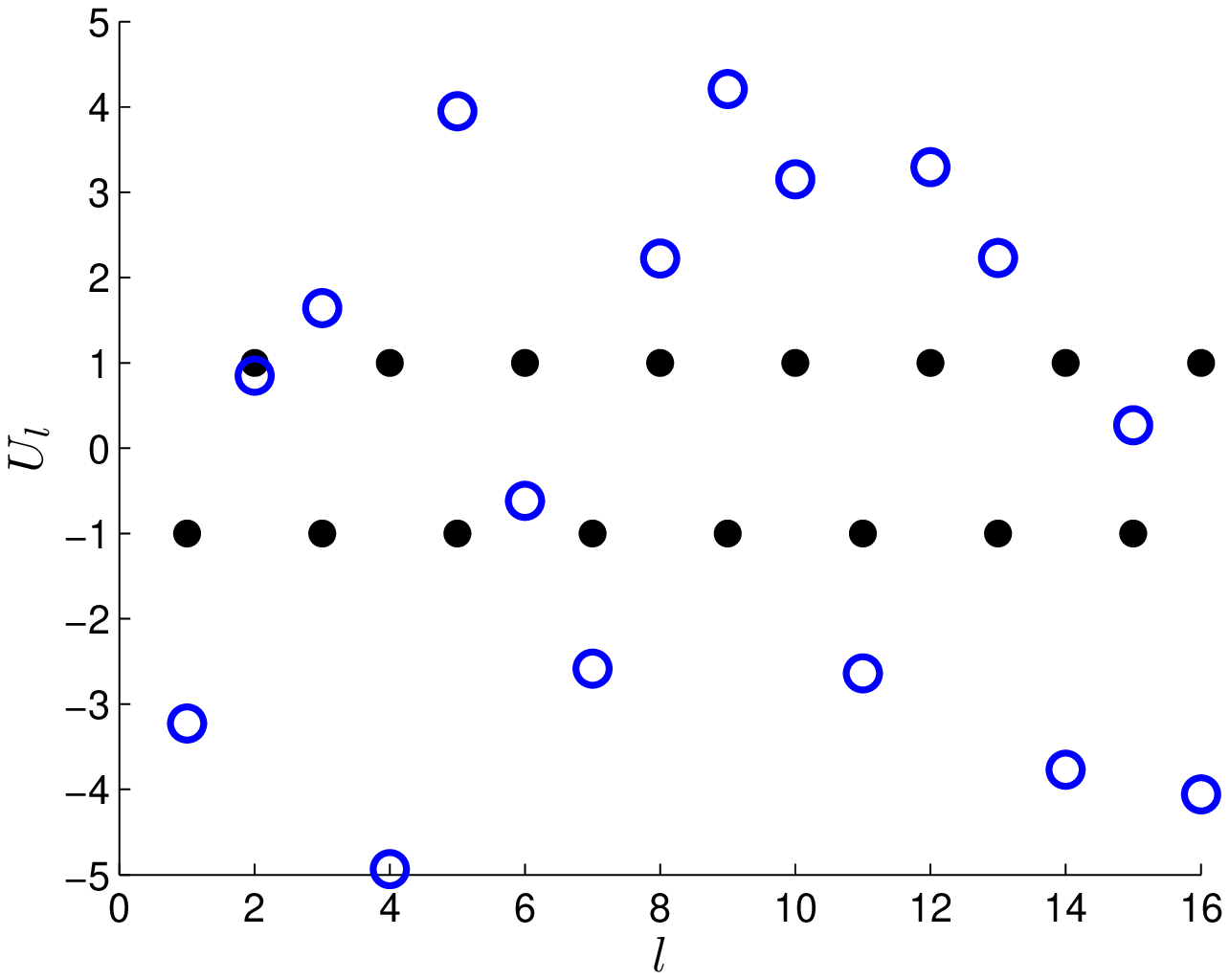}}
  \subfigure[]{\includegraphics[width=0.50\textwidth]{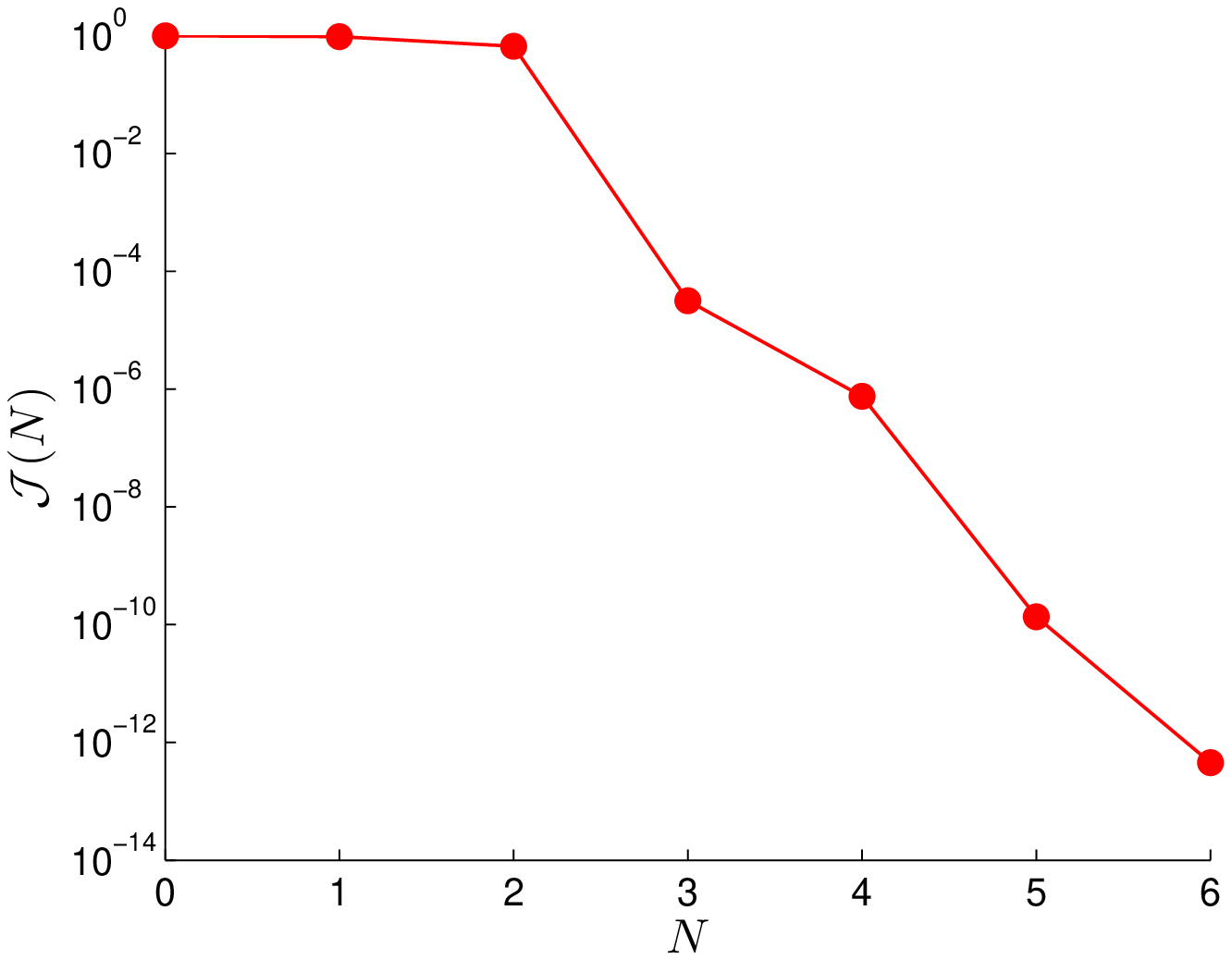}}}
  \caption{(a) Empty blue circles show optimal solution $\hat U$ reconstructed from the initial guess $U_{l,ini}$
    (filled black circles) provided in Table~\ref{tab:curr_model_C}.
    (b) Cost functional $\mathcal{I}(N)$ as a function of optimization iteration $N$ in solving the EIT
    Problem to find optimal solution $(\hat u(x), \hat U)$.}
  \label{fig:eit_res_forward}
  \end{center}
\end{figure}

Numerical result of {\it Stage 1} is demonstrated in a Figure~\ref{fig:eit_res_forward}. 
Electrical currents  $(I_l)_{l=1}^{16}$ specified in
Table~\ref{tab:curr_model_C} are injected through 16 electrodes $E_l$, $l = 1, \ldots, 16$,
and electrical conductivity field $\sigma(x)$ is assumed known, i.e.~$\sigma(x) = \sigma_{true}(x)$.
Figure~\ref{fig:eit_res_forward}(a) shows the optimal solution for control $U$ (empty blue circles)
reconstructed from the initial guess $U_{l,ini}$ (filled black circles) provided in Table~\ref{tab:curr_model_C}.
Fast convergence in 6 iterations as seen in Figure~\ref{fig:eit_res_forward}(b) confirms well-posedness of
the EIT Problem and also uniqueness of the global solution $\hat U$ of the convex Problem $\mathcal{I}$ (see Remark~\ref{convex}).
%By finding optimal control $(\hat u(x), \hat U)$ in Problem $\mathcal{I}$,
%we aim reconstructing boundary voltage vector $\hat U$ which will mimic the measurements
%to be used further in the next step of solving Problem $\mathcal{J}$ or the variational formulation of the Inverse EIT Problem.
%

{\it Stage 2.} Solve Problem $\mathcal{J}$ with limited data $I=(I_l)_{l=1}^{16}$ by the gradient descent
method described in Section~\ref{sec:Gradient Method in Banach Space}, Algorithm~\ref{alg:param_opt} and to recover
optimal control $\left( \sigma_{true}, U_{true} \right)$.

\begin{figure}[!htb]
  \begin{center}
  \mbox{
  \subfigure[]{\includegraphics[width=0.50\textwidth]{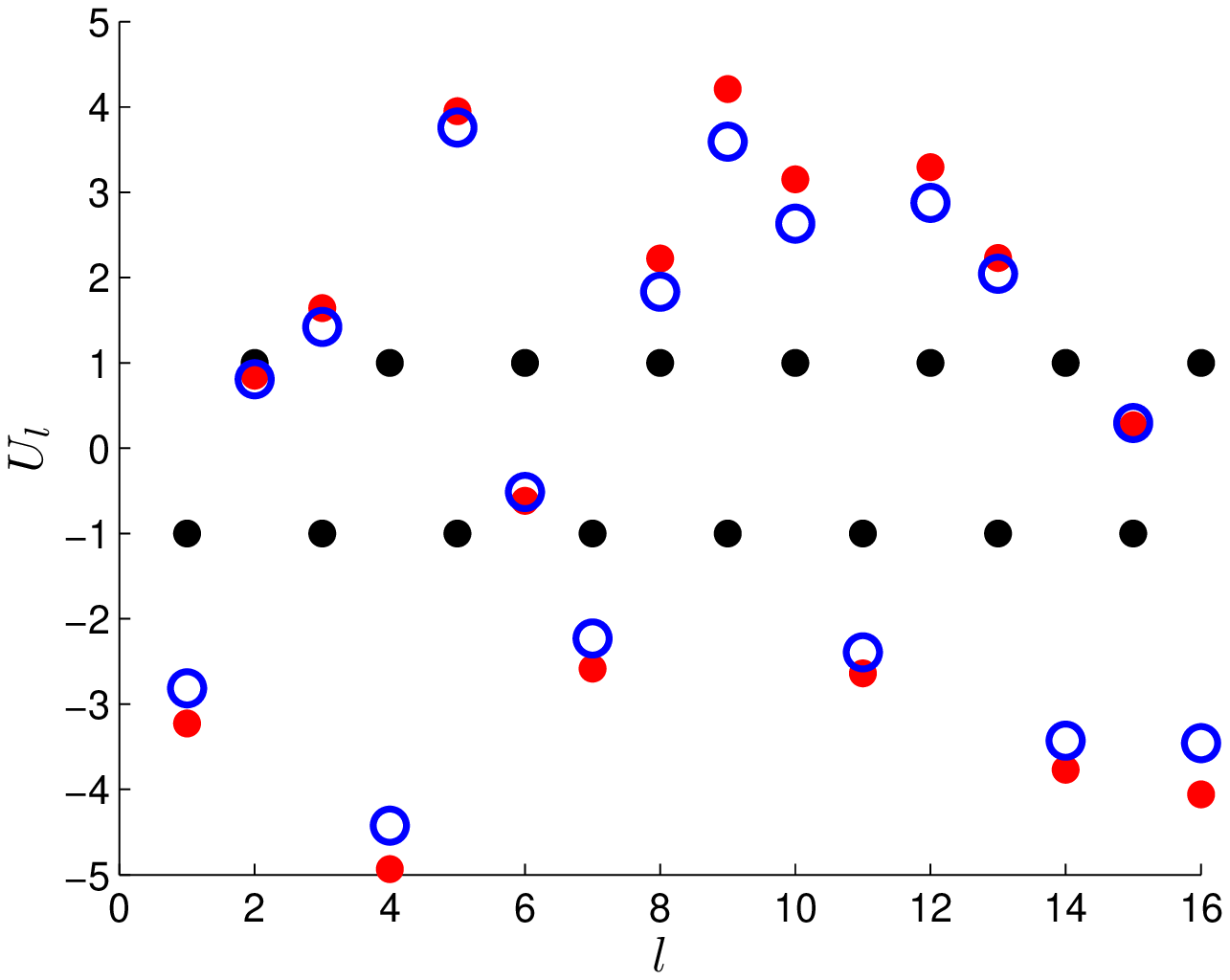}}
  \subfigure[]{\includegraphics[width=0.50\textwidth]{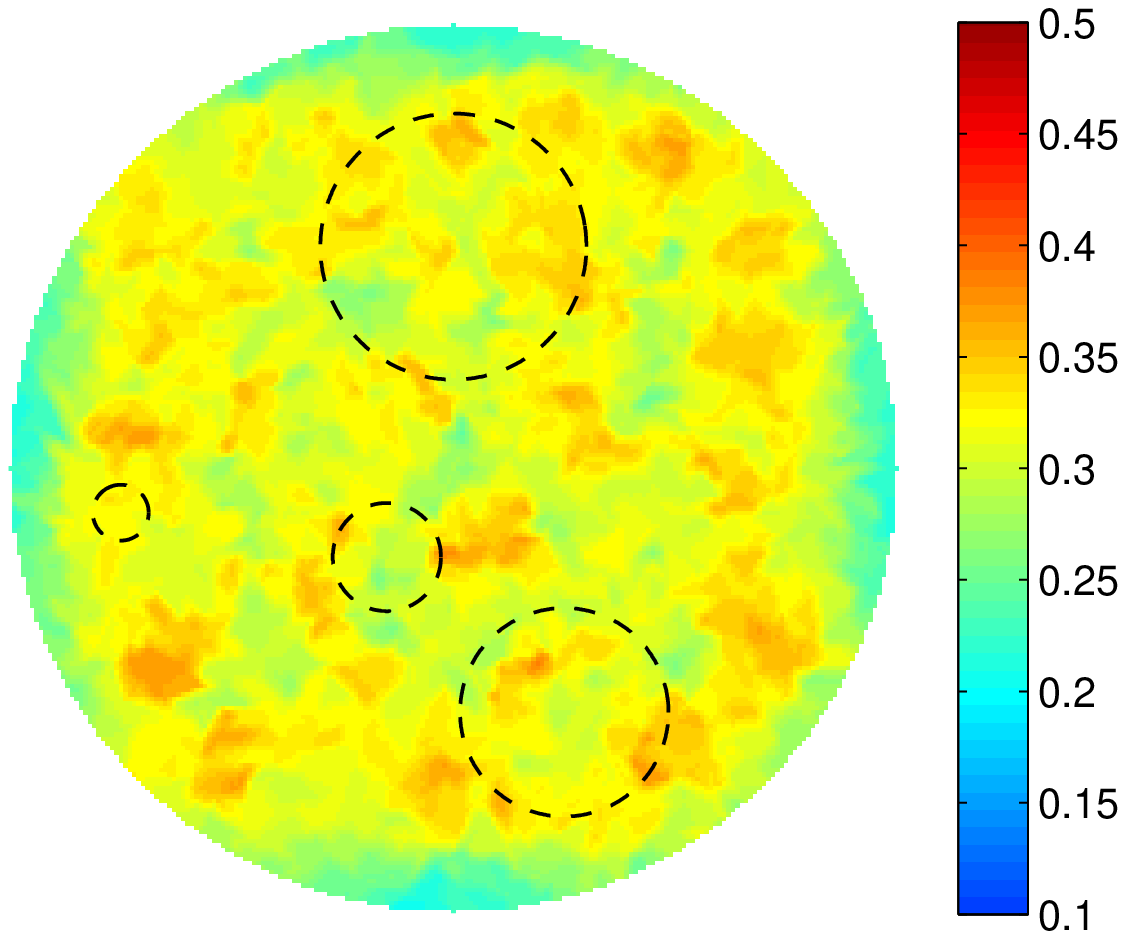}}}
  \caption{(a) Empty blue circles show optimal solution $\hat U$ reconstructed from the initial guess $U_{l,ini}$
    (filled black circles) provided in Table~\ref{tab:curr_model_C}. Filled red circles represent
    actual experimental data $U^*$ (also blue circles in Figure~\ref{fig:eit_res_forward}(a)).
    (b) Reconstructed electrical conductivity field $\hat \sigma(x)$. Dashed circles represent the location
    of four cancer-affected regions taken from known $\sigma_{true}$.
    %Shown optimal solution $(\hat \sigma(x), \hat U)$ is obtained by solving the EIT Inverse Problem without
    %applying additional data acquired through rotating boundary voltages $U_l$.
    }
  \label{fig:eit_res_inv_ell_0_eps_0_noRot}
  \end{center}
\end{figure}

\begin{figure}[!htb]
  \begin{center}
  \mbox{
  \subfigure[]{\includegraphics[width=0.50\textwidth]{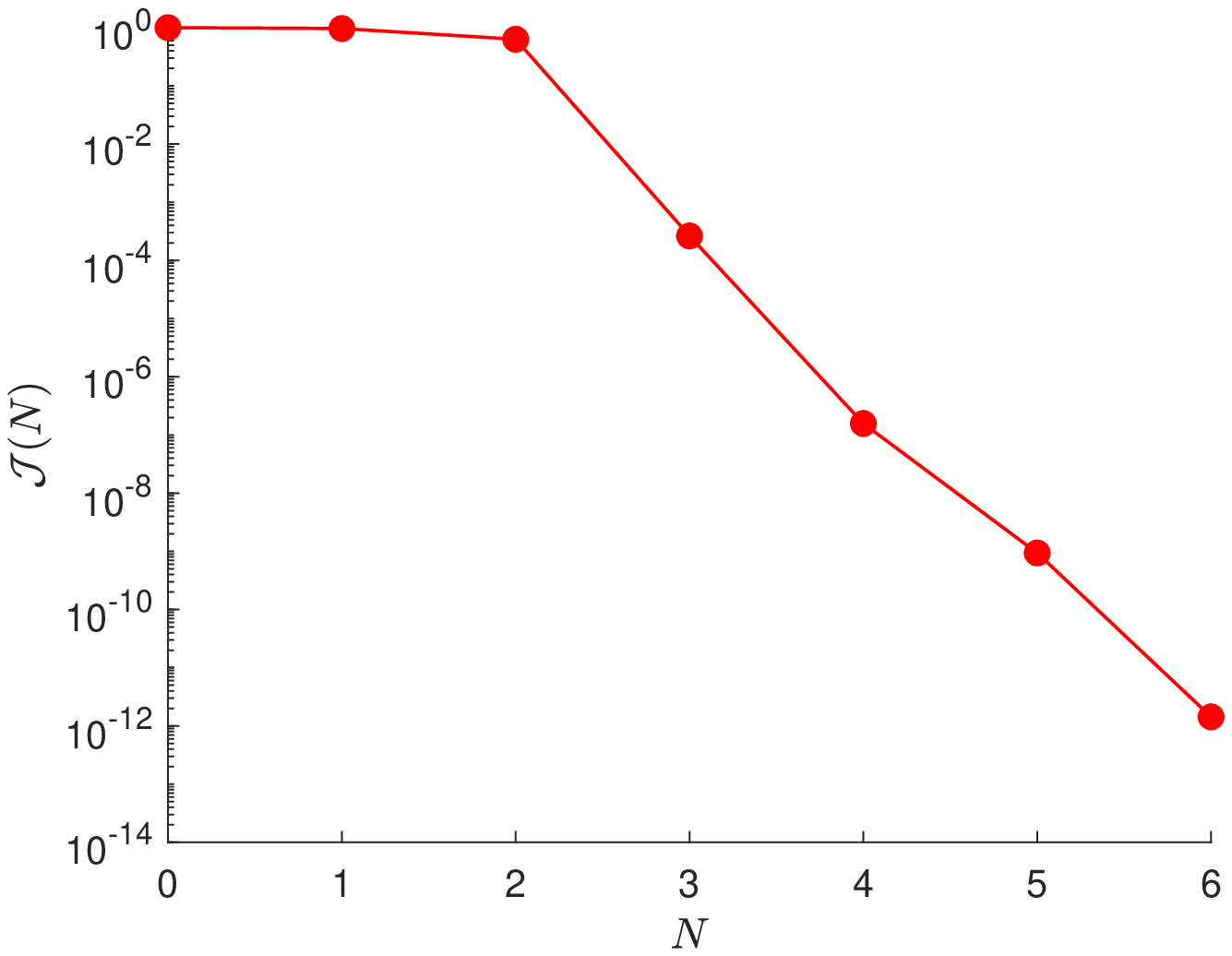}}
  \subfigure[]{\includegraphics[width=0.50\textwidth]{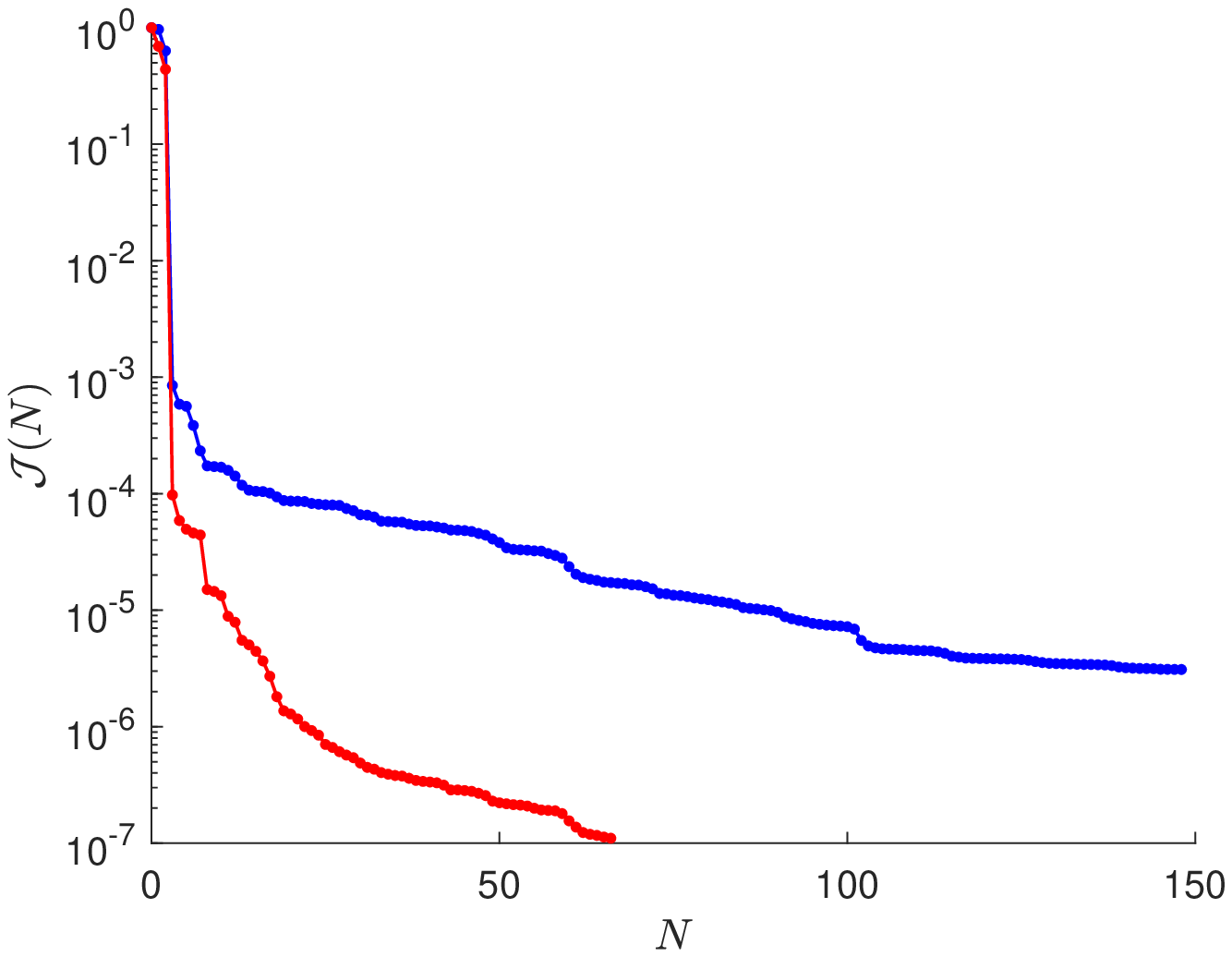}}}
  \caption{Cost functional $\mJ(N)$ as a function of optimization iteration $N$ in solving the EIT Inverse
    Problem to find optimal solution $(\hat \sigma(x), \hat U)$ (a) without and (b) with applying
    additional data acquired through rotating boundary voltages $U_l$. Convergence in (b) is compared for
    two cases: (blue dots) without regularization, and (red dots) when applying regularization with parameter
    $\beta^* = 0.3162$.}
  \label{fig:eit_res_grad_compare}
  \end{center}
\end{figure}

Numerical result of {\it Stage 2} without regularization ($\beta=0$) is demonstrated in a Figure~\ref{fig:eit_res_inv_ell_0_eps_0_noRot}. Furthermore, in all subsequent Figures, we mark the
location of four cancer-affected regions from known $\sigma_{true}$ by dashed circles. As seen in Figure~\ref{fig:eit_res_inv_ell_0_eps_0_noRot}(b),
the electrical conductivity field $\sigma(x)$ is reconstructed poorly without any signature to identify
spots with cancer-affected tissues.  %This simplifies our
%efforts to conclude on how properly these regions are trapped while solving the EIT Inverse Problems via Problem $\mathcal{J}$.
Fast convergence with respect to functional in just 6 iterations is demonstrated in Figure~\ref{fig:eit_res_grad_compare}(a).
However, there is no convergence with respect to all control parameters as shown in Figure~\ref{fig:eit_res_inv_ell_0_eps_0_noRot}(a,b).
Although the $U$-component deviates slightly
from actual experimental data $U^*$ (filled red circles), the optimal solution $\hat \sigma(x)$ obtained for the
$\sigma$-component is significantly different from the true solution $\sigma_{true}$. This is a consequence of the ill-posedness
of the inverse EIT problem due to non-uniqueness of the solution.

{\it Stage 3.} To increase the size of input data we apply the same set of boundary voltages $U^*_l$ to
different electrodes $E_l$ using a ``rotation scheme'', i.e. we denote $U^1=U^*, I^1=I$ and consider 15 new permutations of boundary voltages
as in \eqref{permutations} applied to electrodes $E_1, E_2, \ldots, E_{16}$ respectively. For each boundary voltage vector $U^j$ we solve elliptic PDE problem
\eqref{eq:forward_1}--\eqref{eq:forward_3} to obtain the distribution of electrical potential $u_j(\cdot)=u(\cdot; U^j)$ over boundary $\partial Q$.
By using ``voltage--to--current'' formula  \eqref{eq:forward_4}, we calculate current pattern $I^j$ associated with $U^j$. Thus, a new set $(I^j)^{16}_{j=1}$
contains 256 input data that could be enough to expect the problem to be well-posed in case a reduced-dimensional
space for control $\sigma$ as described in Section~\ref{sec:PCA}. Practical analogy of this step is implementation of the
``voltage--to--current'' procedure: by injecting 15 new sets of voltages $U^j, j=2,...,16$ from \eqref{permutations} on
the electrodes $E_l, \ l= 1, \cdots, 16$, take the measurement of the currents $I^j=(I^j_1,\cdots,I^j_{16})$.
Then we solve Problem $\mathcal{K}$ with extended data set by the gradient descent
method described in Section~\ref{sec:Gradient Method in Banach Space}, Algorithm~\ref{alg:param_opt} and to recover
optimal control $\left( \sigma_{true}, U_{true} \right)$.

\begin{figure}[!htb]
  \begin{center}
  \mbox{
  \subfigure[]{\includegraphics[width=0.50\textwidth]{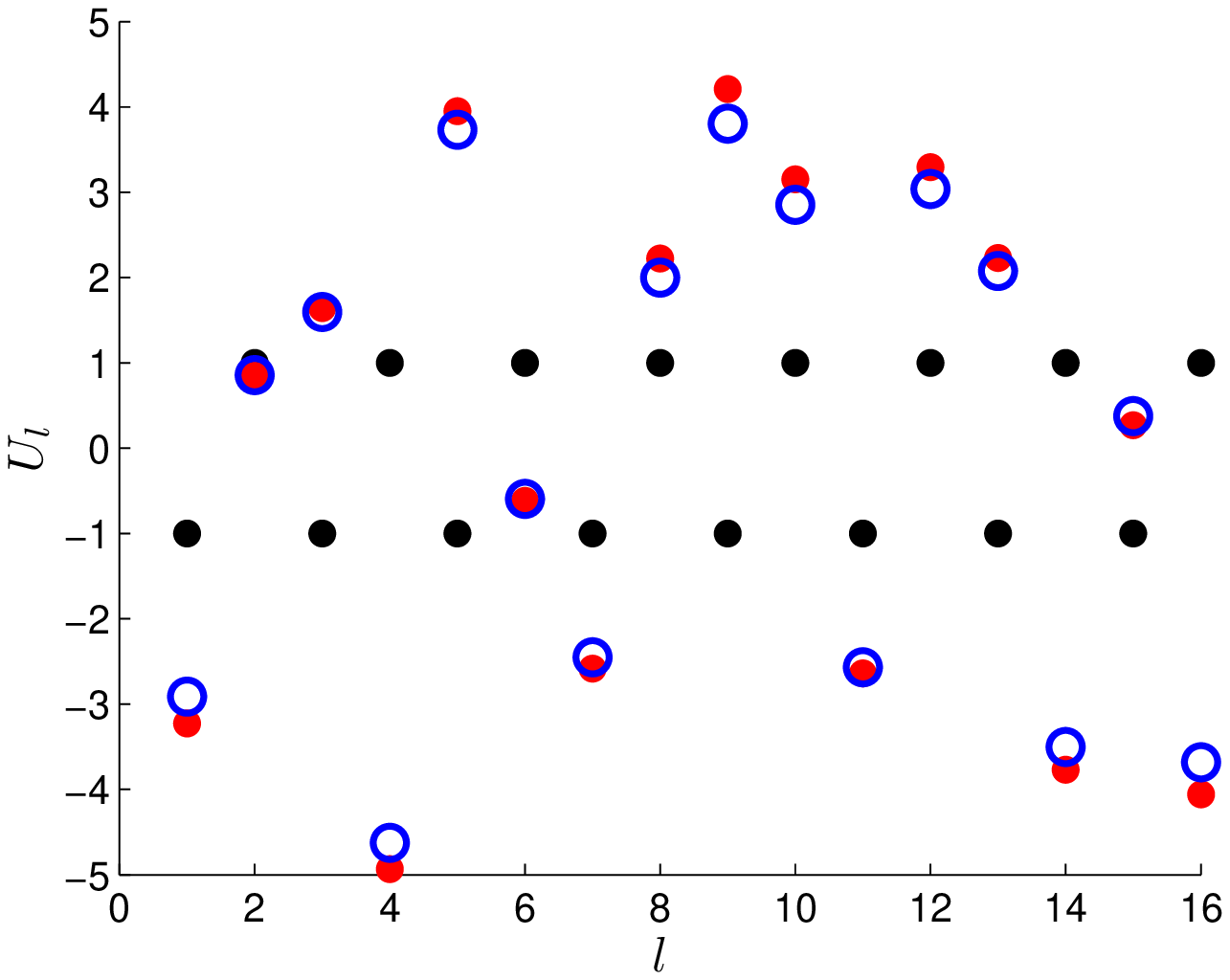}}
  \subfigure[]{\includegraphics[width=0.50\textwidth]{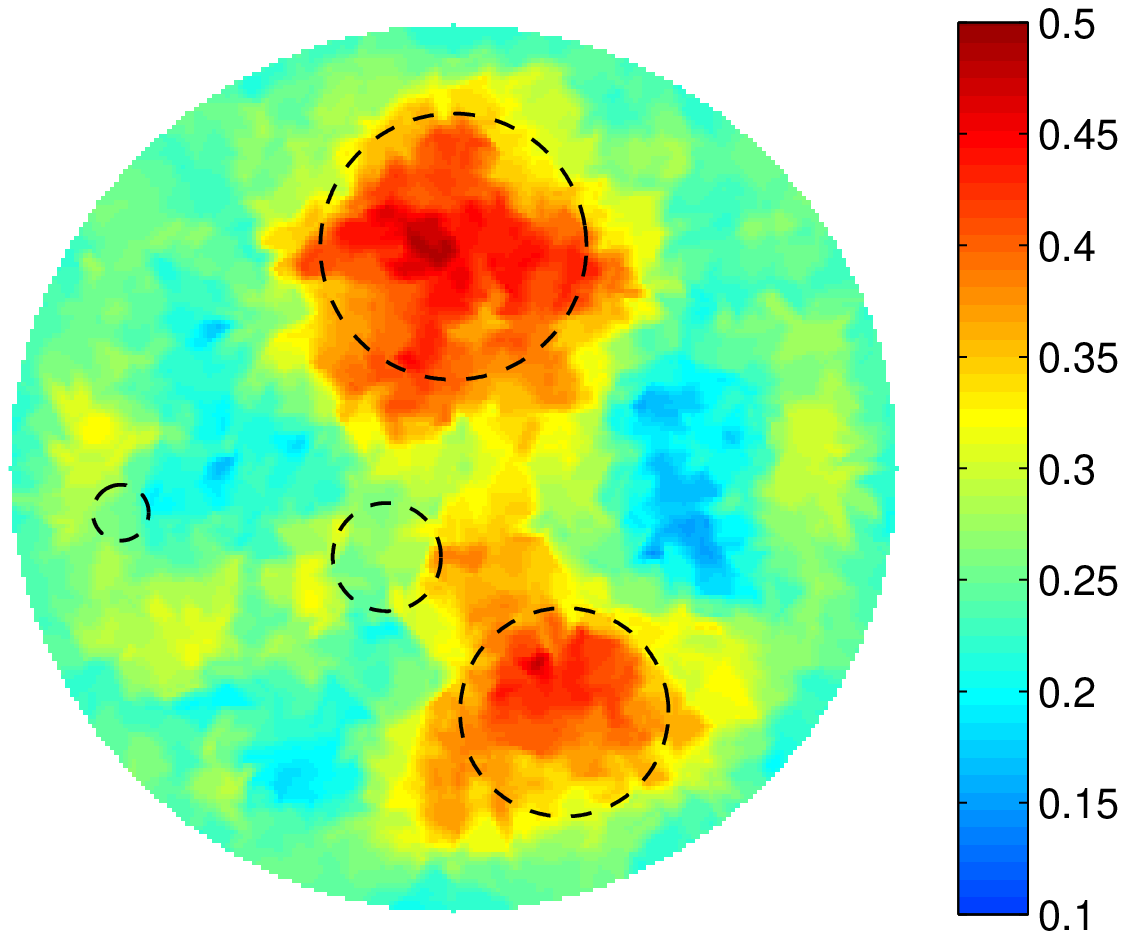}}}
  \caption{(a) Empty blue circles show optimal solution $\hat U$ reconstructed from the initial guess $U_{l,ini}$
    (filled black circles) provided in Table~\ref{tab:curr_model_C}. Filled red circles represent
    actual experimental data $U^*$ (also blue circles in Figure~\ref{fig:eit_res_forward}(a)).
    (b) Reconstructed electrical conductivity field $\hat \sigma(x)$. Dashed circles represent the location
    of four cancer-affected regions taken from known $\sigma_{true}$.
    %Shown optimal solution $(\hat \sigma(x), \hat U)$ is obtained by solving the EIT Inverse Problem with
    %applying additional data acquired through rotating boundary voltages $U_l$.
    }
  \label{fig:eit_res_inv_ell_0_eps_0}
  \end{center}
\end{figure}

Numerical result of {\it Stage 3} without regularization ($\beta=0$) is demonstrated in a Figure~\ref{fig:eit_res_inv_ell_0_eps_0}. Contrary to previous results, 
the electrical conductivity field $\sigma(x)$ is reconstructed
much better matching the two biggest spots while not perfectly capturing their shapes. Reconstruction result for boundary
voltage $U$ is also improved.
%In the next two Figures we are going to compare these optimization results with those obtained
%by employing two numerical approaches discussed in Sections~\ref{sec:tuning_pca}--\ref{sec:bc_smoothing} and
%discuss any changes in performance.
%

\begin{figure}[!htb]
  \begin{center}
  \mbox{
  \subfigure[]{\includegraphics[width=0.50\textwidth]{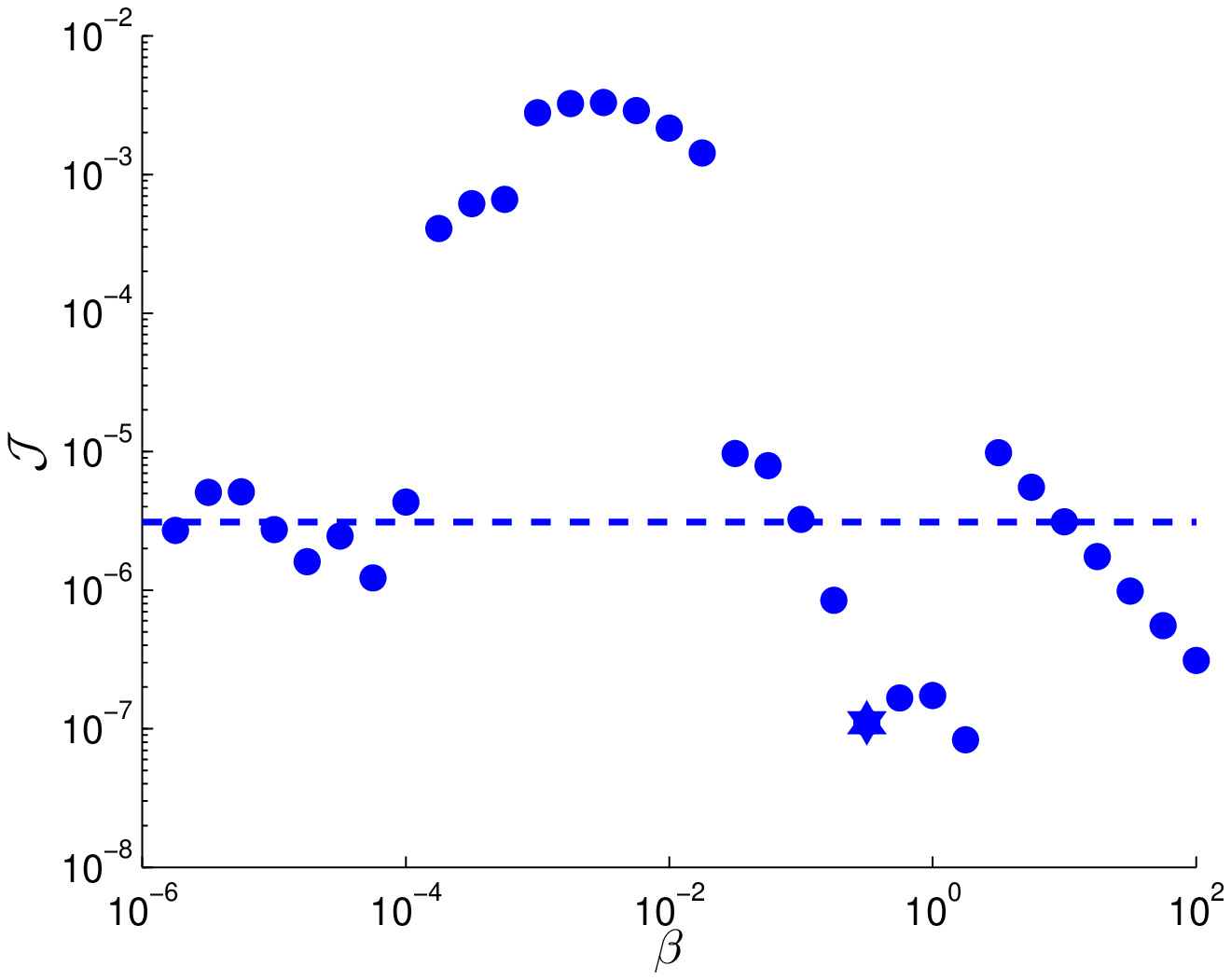}}
  \subfigure[]{\includegraphics[width=0.50\textwidth]{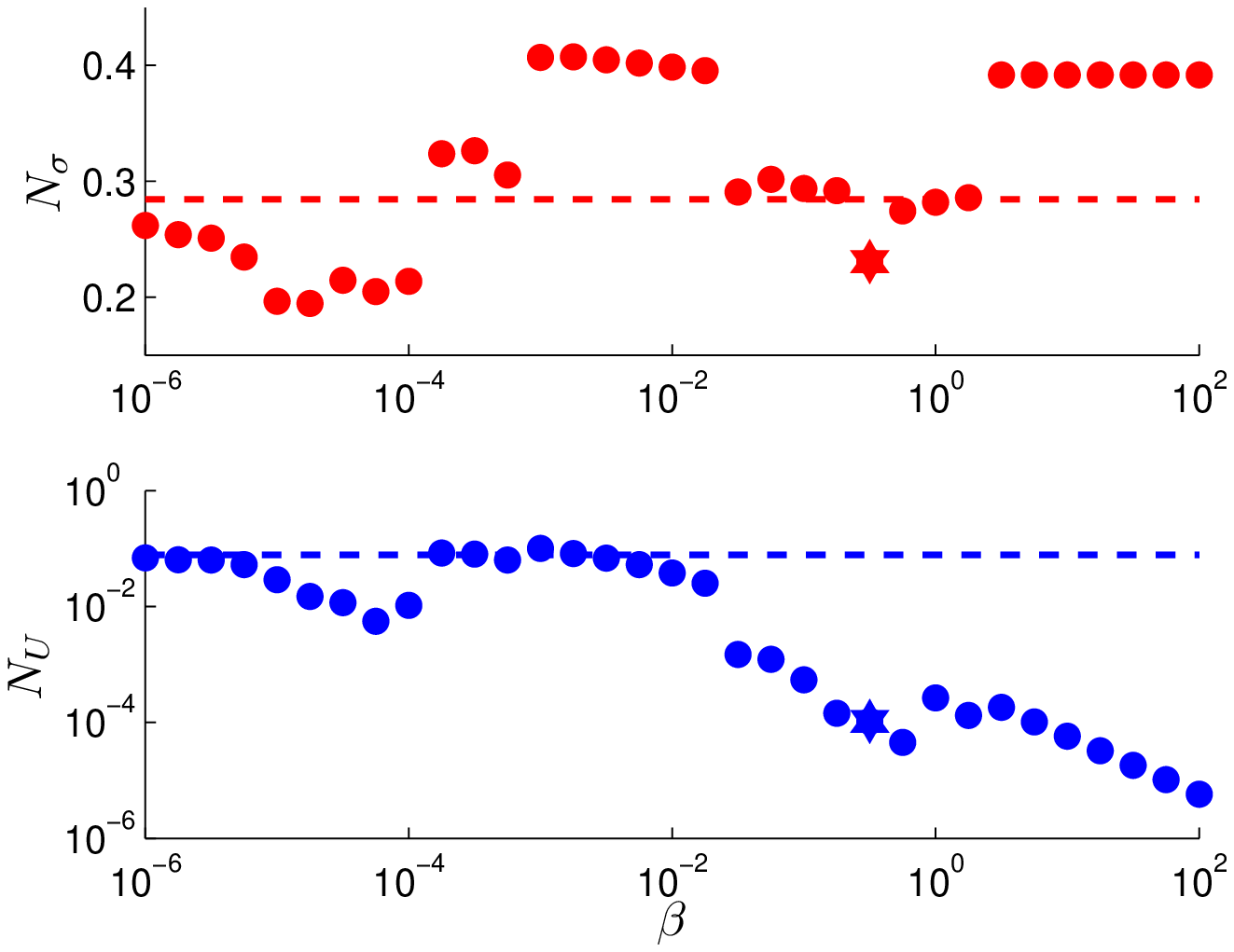}}}
  \caption{(a) Cost functional $\mathcal{K}$ values and (b) solution norms
    $N_{\sigma} = \frac{\| \sigma - \sigma_{true}\|_{L_2}}{\| \sigma_{true}\|_{L_2}}$ and
    $N_U = \frac{| U - U^*|}{| U^*|}$ evaluated at termination (dots) for different
    values of regularization parameter $\beta$ in \eqref{eq:cost_functional} and (dashed lines) when $\beta=0$.
    The best results obtained at $\beta^* = 0.3162$ are shown by hexagons.}
  \label{fig:reg_beta_opt}
  \end{center}
\end{figure}

\begin{figure}[!htb]
  \begin{center}
  \mbox{
  \subfigure[]{\includegraphics[width=0.50\textwidth]{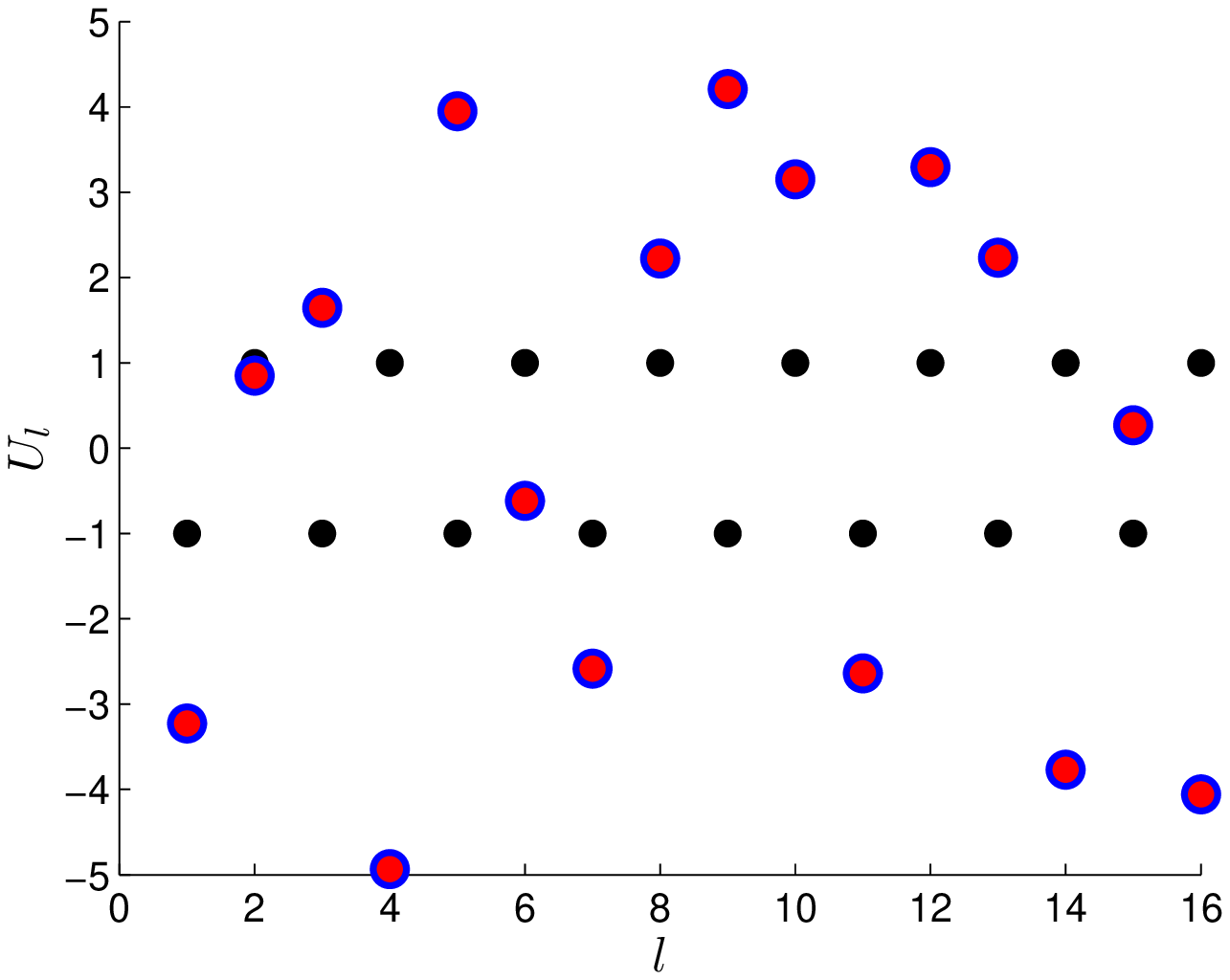}}
  \subfigure[]{\includegraphics[width=0.50\textwidth]{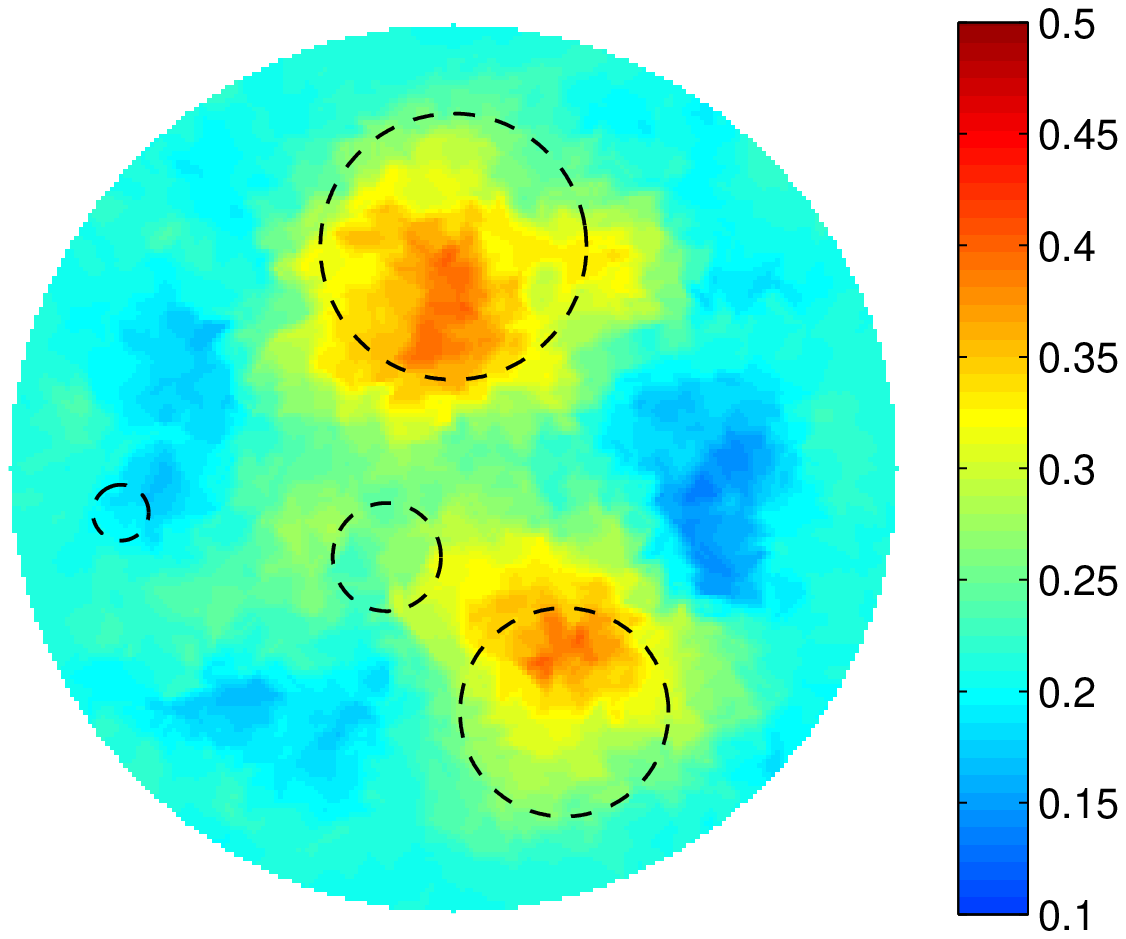}}}
  \caption{(a) Empty blue circles show optimal solution $\hat U$ reconstructed from the initial guess $U_{l,ini}$
    (filled black circles) provided in Table~\ref{tab:curr_model_C}. Filled red circles represent
    actual experimental data $U^*$ (also blue circles in Figure~\ref{fig:eit_res_forward}(a)).
    (b) Reconstructed electrical conductivity field $\hat \sigma(x)$. Dashed circles represent the location
    of four cancer-affected regions taken from known $\sigma_{true}$.
    Optimal solution $(\hat \sigma(x), \hat U)$ is obtained by solving the Problem $\mathcal{K}$ with
    regularization parameter $\beta^* = 0.3162$ in \eqref{eq:cost_functional_rotation}.}
  \label{fig:eit_res_inv_beta_opt}
  \end{center}
\end{figure}

Finally, we evaluate the effect of adding regularization term ($\beta > 0$) in the cost functional \eqref{eq:cost_functional_rotation}.
The outcomes with respect to different values of regularization parameter $\beta$ (blue dots) are shown in
Figure~\ref{fig:reg_beta_opt}(a). The dashed line represents the result of optimization with $\beta=0$. Numerical results demonstrate that 
small values of $\beta$ (roughly when $\beta < 10^{-4}$) have no significant effect
towards decreasing the values of the cost functional $\mathcal{K}$. Significant improvement at
different scales is observed when $\beta > 10^{-1}$.  To identify optimal value for $\beta$,
we examine additionally $\sigma$ and $U$ solution norms
$N_{\sigma} = \frac{\| \sigma - \sigma_{true}\|_{L_2}}{\| \sigma_{true}\|_{L_2}}$ and
$N_U = \frac{| U - U^*|}{| U^*|}$ presented in Figure~\ref{fig:reg_beta_opt}(b). Based on the numerical results, we pick up the value
(shown by hexagons)  $\beta^* = 0.3162$ as the best value in terms of
improvement of solutions simultaneously with respect to both controls $\sigma$ and $U$. Figure~\ref{fig:eit_res_inv_beta_opt}
shows optimal solution $(\hat \sigma(x), \hat U)$ obtained by choosing  $\beta^* = 0.3162$. Overall optimization performance in
the last case is also enhanced by much faster convergence.
Figure~\ref{fig:eit_res_grad_compare}(b) provides the comparison for convergence results obtained for two different cases,
namely without regularization (blue dots), and with regularization with parameter $\beta^* = 0.3162$ (red dots).

\section{Conclusions}
\label{sec:conclusions}
This paper analyzes the inverse EIT problem on recovering electrical conductivity tensor and potential in the body based on the measurement of the boundary voltages on the electrodes for a given electrode current. The inverse EIT problem presents an effective mathematical model of breast cancer detection based on the experimental fact that the electrical conductivity of malignant tumors of the breast may significantly differ from conductivity of the surrounding normal tissue. We analyze the inverse EIT problem in a PDE constrained optimal control framework in Besov space, where the electrical conductivity tensor and boundary voltages are control parameters, and the cost functional is the norm declinations of the boundary electrode current from the given current pattern and boundary electrode voltages from the measurements. The state vector is a solution of the second order elliptic PDE in divergence form with bounded measurable coefficients under mixed Neumann/Robin type boundary condition. The following are the main results of the paper:
\begin{itemize}
\item In contrast with the current state of the field, the inverse EIT problem is investigated with unknown electrical conductivity tensor, which is essential in understanding 
and detecting the highly anisotropic distribution of cancerous tumors in breast tissue.
\item To address the highly ill-posed nature of the inverse EIT problem, we develop a "variational formulation with additional data" which is well adapted to clinical situation when additional ``voltage--to--current" measurements significantly increase the size of the input data while keeping the size of the unknown parameters fixed.
\item  Existence of the optimal control and Fr\'echet differentiability in the Besov space setting is proved. The formula for the Fr\'echet gradient and optimality condition is derived. 
Effective numerical method based on the projective gradient method in Besov spaces is developed. 
\item Extensive numerical analysis is pursued in the 2D case through implementation of the projective gradient method, re-parameterization via PCA, and Tikhonov regularization in a carefully 
constructed model example which adequately represents the diagnosis of breast cancer in reality. Numerical analysis demonstrates accurate reconstruction of the electrical conductivity
function of the body in the frame of the model based on "variational formulation with additional data".
\end{itemize}
%\section*{Acknowledgements}

% references
%\newpage
%\nocite{*}
\bibliographystyle{agsm}
\bibliography{biblio_Abdulla}

\section*{Appendices}
\label{sec:appendix}

\addcontentsline{toc}{chapter}{Appendices}
\renewcommand{\thesection} {\Alph{section}}
\setcounter{section}{0}
\numberwithin{equation}{section}
\numberwithin{theorem}{section}
\numberwithin{figure}{section}

\section{Validation of Gradients}
\label{sec:kappa_test}

In this section we present results demonstrating the consistency of cost functional gradients
$\nabla_{\sigma} \mJ$, $\nabla_{\xi} \mJ$ and $\nabla_U \mJ$ obtained with the approach described
in Section~\ref{sec:Main_Results} and Algorithm~\ref{alg:param_opt}. As the sensitivity of cost
functionals $\mJ(\sigma, U)$ and $\mJ(\xi, U)$ with respect to controls may vary significantly for
different contributions of $\sigma$, $\xi$ and $U$, it is reasonable to perform testing separately
for different parts of the gradients, namely $\nabla_{\sigma} \mJ(\sigma,U)$,
$\nabla_{\xi} \mJ(\xi,U)$ and $\nabla_U \mJ(\sigma,U)$.

\begin{figure}[htb!]
  \begin{center}
  \mbox{
  \subfigure[]{\includegraphics[width=0.5\textwidth]{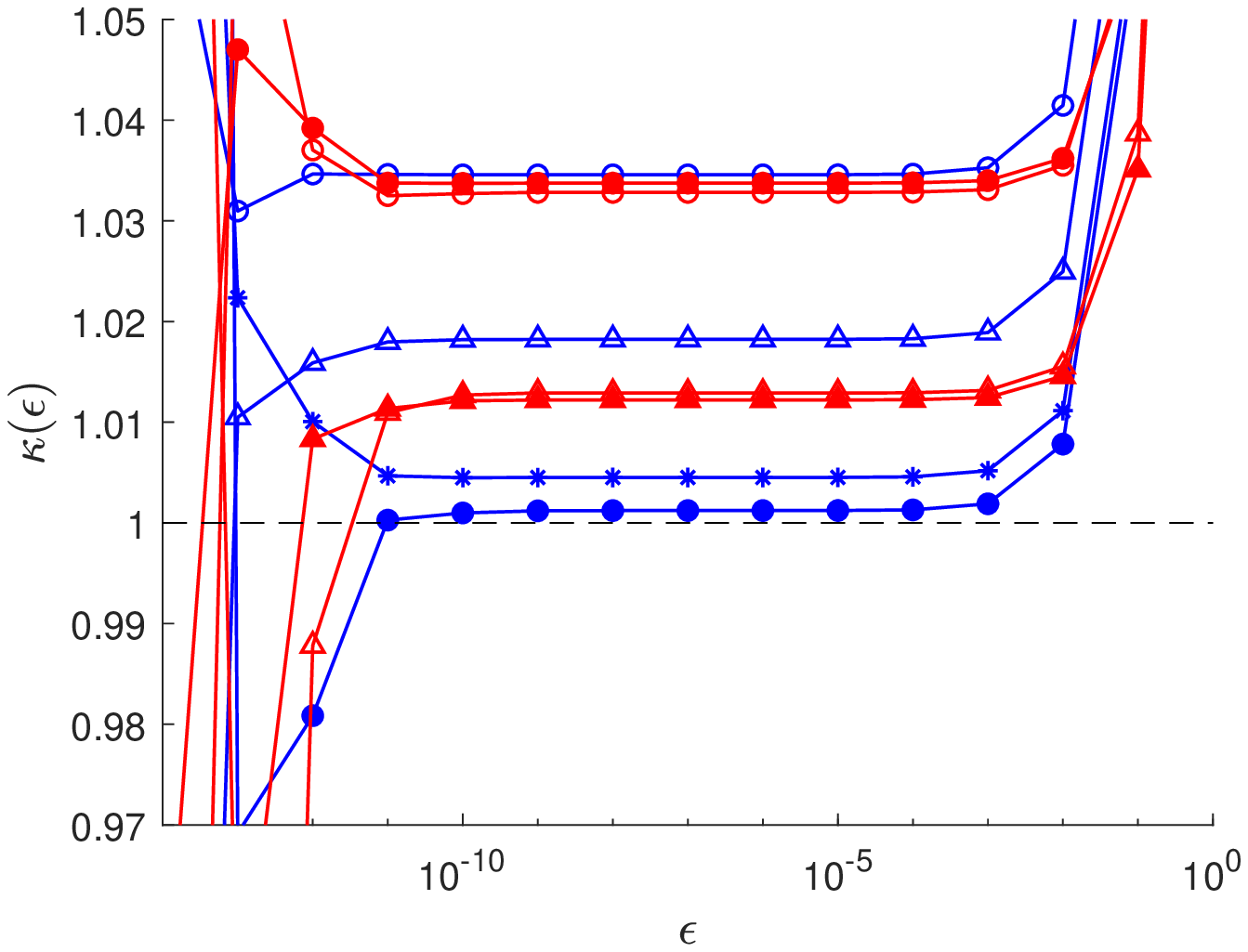}}
  \subfigure[]{\includegraphics[width=0.5\textwidth]{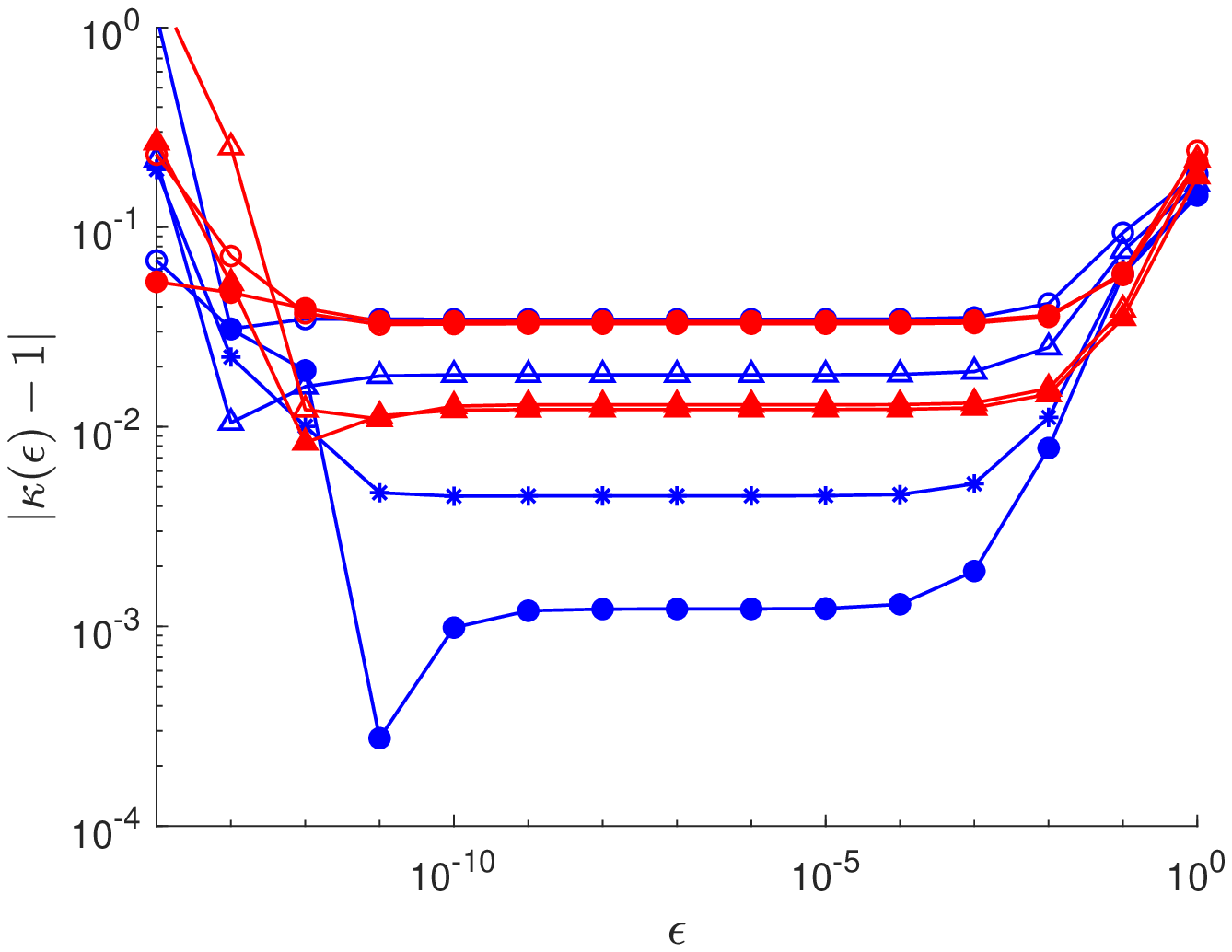}}}
  \end{center}
  \caption{The behavior of (a) $\kappa (\epsilon)$ and (b) $\log_{10} |\kappa(\epsilon) -1 |$ as a function of $\epsilon$
    for both controls (blue) $\sigma(x)$ and (red) $\xi$. Number of vertices over boundary $\partial Q$ used in spatial
    discretization of domain $Q$ for computing $\nabla_{\sigma} \mJ(\sigma,U)$ are
    (open circles) $n_v = 48$ (724 mesh $\mM$ triangular finite elements), (triangles) $n_v = 96$ (1996 elements),
    (asterisks) $n_v = 176$ (7800 elements), (filled circles) $n_v = 432$ (47938 elements).
    Same meshes (red circles) $\mM(48)$ and (red triangles) $\mM(96)$ are used for computing $\nabla_{\xi} \mJ(\xi,U)$
    with number of principal components (open circles/triangles) $N_{\xi} = 20$ and (filled circles/triangles) $N_{\xi} = 74$.}
  \label{fig:kappa_test}
\end{figure}

First, we explore the results obtained for controls representing electrical conductivity $\sigma$
before and after projecting the gradients onto the reduced-dimensional $\xi$-space as described
in Section~\ref{sec:PCA}. Figure~\ref{fig:kappa_test} shows the results of a diagnostic test commonly
employed to verify correctness of cost functional gradients (see, e.g., \cite{Bukshtynov11,Bukshtynov13})
computed for our computational model detailed in Section~\ref{sec:model}. Testing
$\nabla_{\sigma} \mJ(\sigma,U)$ consists in computing the Fr\'echet differential
$d\mJ(\sigma,U; \delta \sigma) = \left\langle{}\J'(\sigma,U),{\delta \sigma}\right\rangle_H$
for some selected variations (perturbations) $\delta \sigma$ in two different ways, namely, using
a finite--difference approximation and using \eqref{eq:Frechet_Derivative} which is based on the
adjoint field, and then examining the ratio of the two quantities, i.e.,
\begin{equation}
  \kappa (\epsilon) = \dfrac{\frac{1}{\epsilon} \left[ \J(\sigma + \epsilon \, \delta \sigma,U) - \J(\sigma,U) \right]}
  { \left\langle{}\J'(\sigma,U),{\delta \sigma}\right\rangle_H}
  \label{eq:kappa}
\end{equation}
for a range of values of $\epsilon$.  If these gradients are computed correctly, then for intermediate values of
$\epsilon$, $\kappa(\epsilon)$ will be close to the unity. Remarkably, this behavior can be observed in
Figure~\ref{fig:kappa_test}(a) over a range of $\epsilon$ spanning about 8-9 orders of magnitude for controls $\sigma$.
Furthermore, we also emphasize that refining mesh $\mM(n_v)$ in discretizing domain $Q$ while solving both state \eqref{eq:forward_1}--\eqref{eq:forward_3} and adjoint \eqref{adj pde}--\eqref{bdry cond of adj pde} PDE problems
yields values of $\kappa (\epsilon)$ closer to the unity. The reason is that in the ``optimize--then--discretize''
paradigm adopted here such refinement of discretization leads to a better approximation of the continuous gradient
as shown in \cite{ProtasBewleyHagen04}. We add that the quantity $\log_{10} | \kappa (\epsilon) - 1 |$ plotted in
Figure~\ref{fig:kappa_test}(b) shows how many significant digits of accuracy are captured in a given gradient evaluation.
As can be expected, the quantity $\kappa(\epsilon)$ deviates from the unity for very small values of $\epsilon$, which
is due to the subtractive cancelation (round--off) errors, and also for large values of $\epsilon$, which is due to
the truncation errors, both of which are well--known effects.

The same test could be easily applied for controls $\xi$ to check the consistency for gradients $\nabla_{\xi} \mJ(\xi,U)$
in the reduced-dimensional $\xi$-space. As seen in Figure~\ref{fig:kappa_test} the same conclusion could be made on the
effect of refining mesh $\mM(n_v)$ in discretizing domain $Q$. We should notice that applying PCA-based re-parameterization
improves the results of this diagnostics. At the same time we conclude that changing number of principal components
$N_{\xi}$ influences the results of the test insignificantly. This effect is easily explained by the fact that only
the first, and thus the biggest, components have sufficient weight and prevails over the rest components in the truncated
tail of the PCA-component sequence.

\begin{figure}[htb!]
  \begin{center}
  \mbox{
  \subfigure[]{\includegraphics[width=0.5\textwidth]{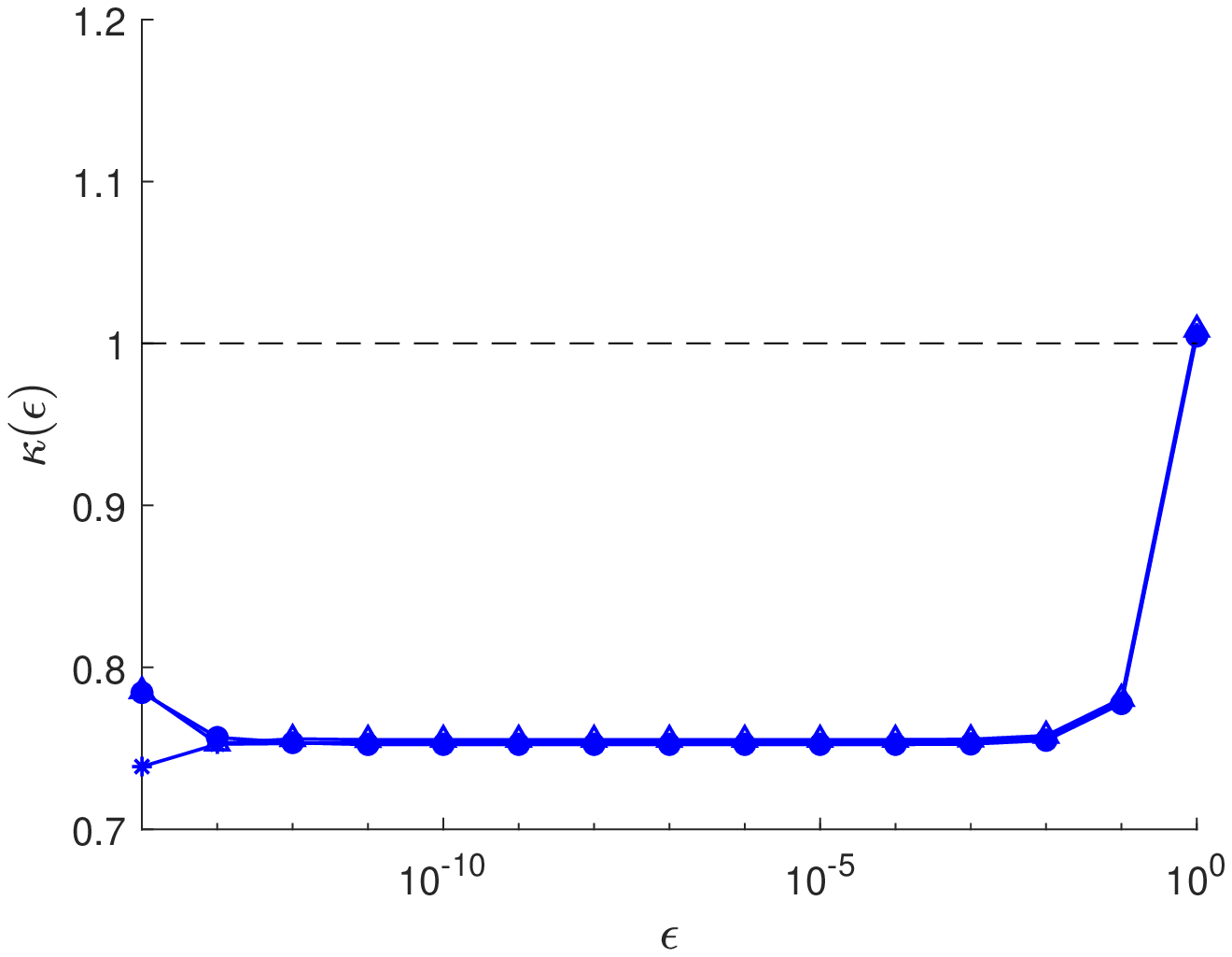}}
  \subfigure[]{\includegraphics[width=0.5\textwidth]{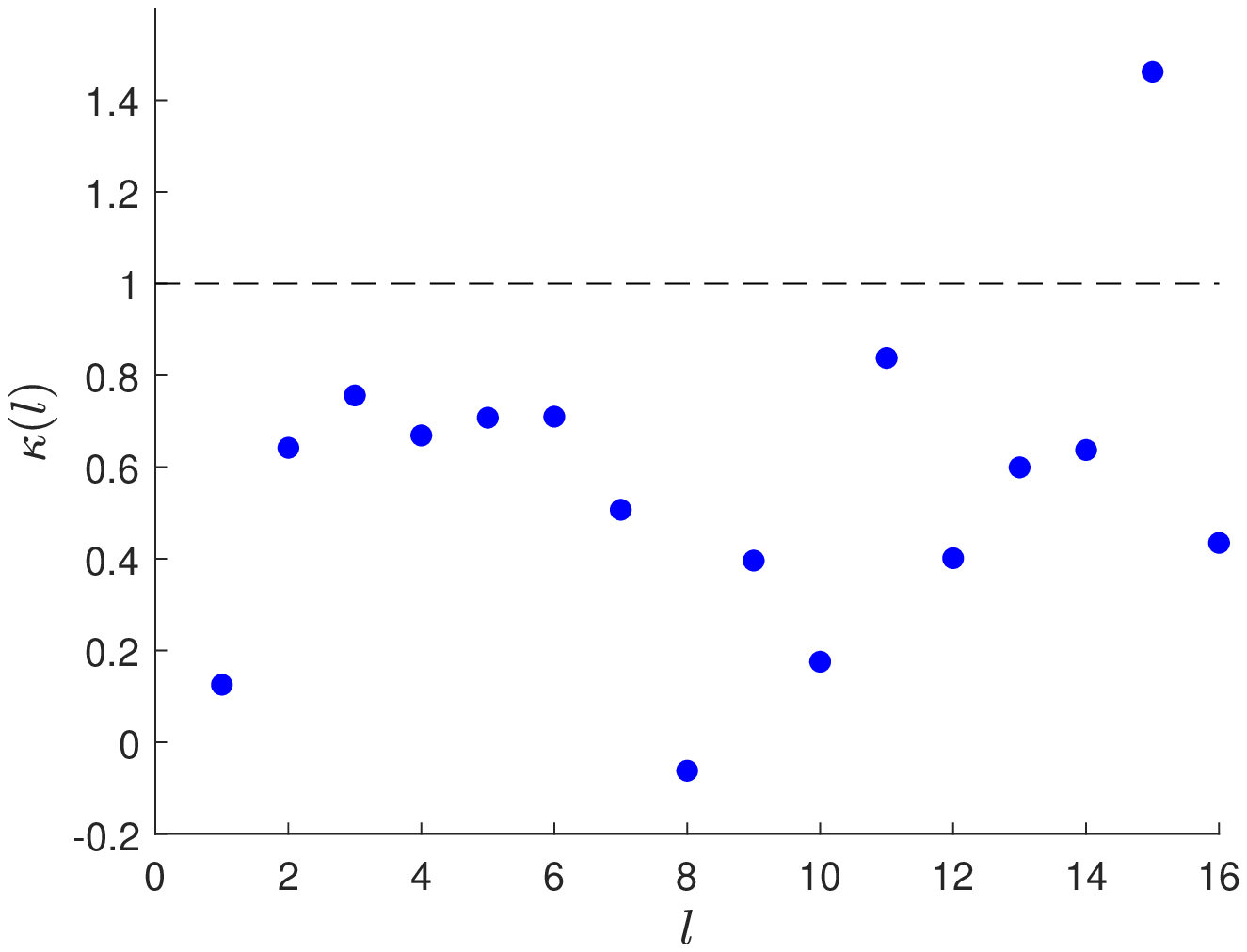}}}
  \end{center}
  \caption{(a) The behavior of $\kappa (\epsilon)$ as a function of $\epsilon$ for control $U$. Number of vertices over
    boundary $\partial Q$ used in spatial discretization of domain $Q$ for computing $\nabla_U \mJ(\sigma,U)$ are
    (triangles) $n_v = 96$ (1996 mesh $\mM$ triangular finite elements), (asterisks) $n_v = 176$ (7800 elements),
    (filled circles) $n_v = 432$ (47938 elements). (b) The behavior of $\kappa(l)$ as a function of electrode number~$l$
    for control $U$ computed by \eqref{eq:kappa_U}.}
  \label{fig:kappa_test_U}
\end{figure}

Second, we could also apply the same testing technique to check the correctness and consistency for gradients
$\nabla_U \mJ(\sigma,U)$ as shown in Figure~\ref{fig:kappa_test_U}(a). Unlike for tests performed for controls
$\sigma$ and $\xi$, gradients $\nabla_U \mJ(\sigma,U)$ are computed correctly but with much larger error demonstrated by
the plateau form of $\kappa (\epsilon)$ which is quite distant from the unity. Remarkably, refining mesh $\mM(n_v)$
in discretizing domain $Q$ while solving both state \eqref{eq:forward_1}--\eqref{eq:forward_3} and
adjoint \eqref{adj pde}--\eqref{bdry cond of adj pde} PDE problems does not change significantly the quality of the
obtained gradients with respect to controls $U$. We explain this by the fact that computing gradients $\nabla_U \mJ(\sigma,U)$
relies mainly on the solution for the potential $u(x)$ obtained on or very close to boundary $\partial Q$ where it looses its
regularity due to discontinuous boundary conditions \eqref{eq:forward_2}--\eqref{eq:forward_3}. As control vector $U$
contains only $m=16$ components, we could also perform our diagnostic test applied individually to every component $U_l$
for fixed (intermediate) value of $\epsilon$
\begin{equation}
  \kappa (l) = \dfrac{\frac{1}{\epsilon} \left[ \J(\sigma, U_l + \epsilon \, \delta U_l) - \J(\sigma,U_l) \right]}
  { \J'(\sigma,U_l) \cdot \delta U_l}.
  \label{eq:kappa_U}
\end{equation}
Figure~\ref{fig:kappa_test_U}(b) represents the results of this modified test which may also be used in analysis
for sensitivity of cost functional $\mJ(\sigma,U)$ to changes in boundary potential $U_l$ at individual electrode $E_l$.

\section{Optimal Size of Reduced-Dimensional $\xi$-space}
\label{sec:tuning_pca}

\begin{figure}[htb!]
  \begin{center}
  \mbox{
  \subfigure[]{\includegraphics[width=0.5\textwidth]{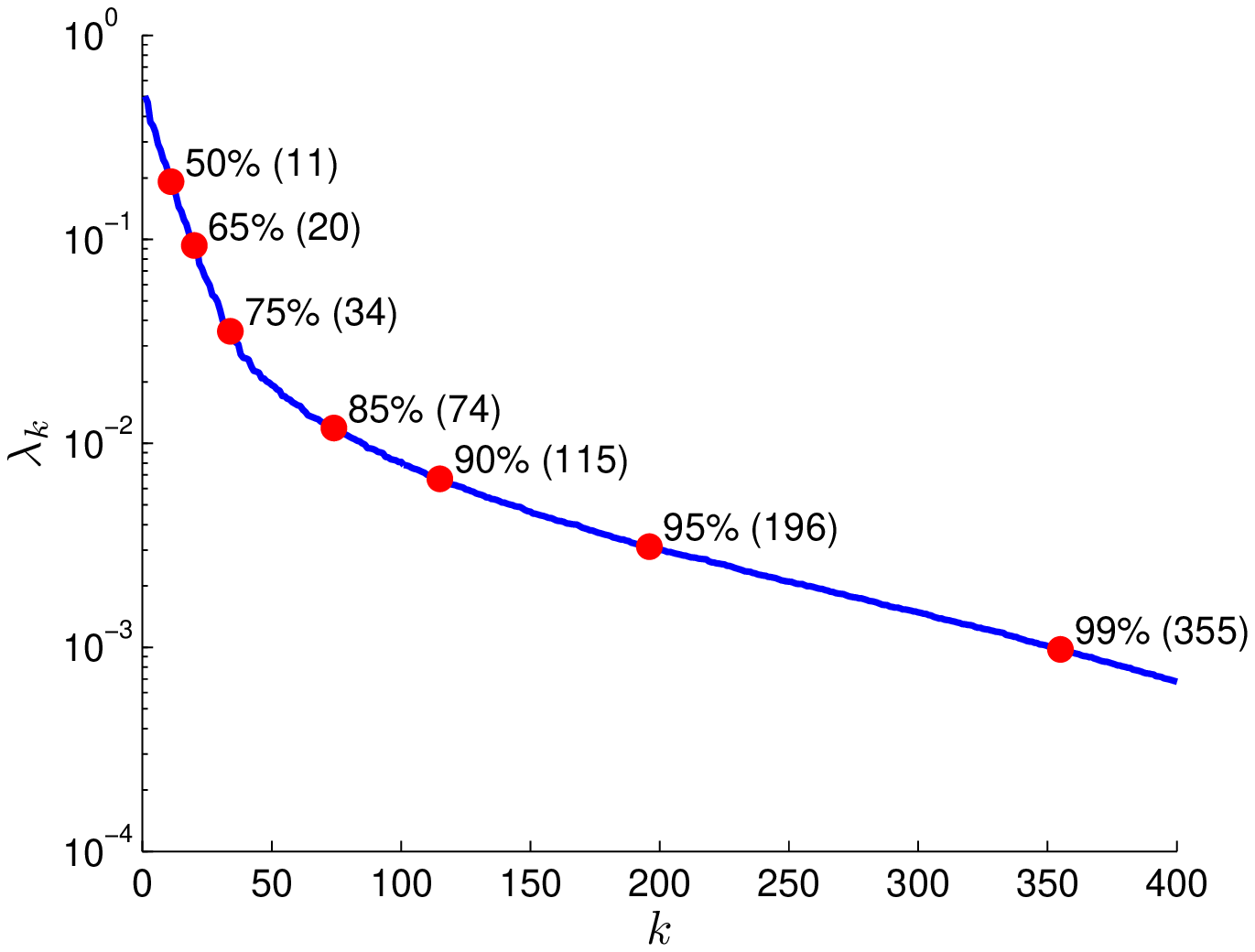}}
  \subfigure[]{\includegraphics[width=0.5\textwidth]{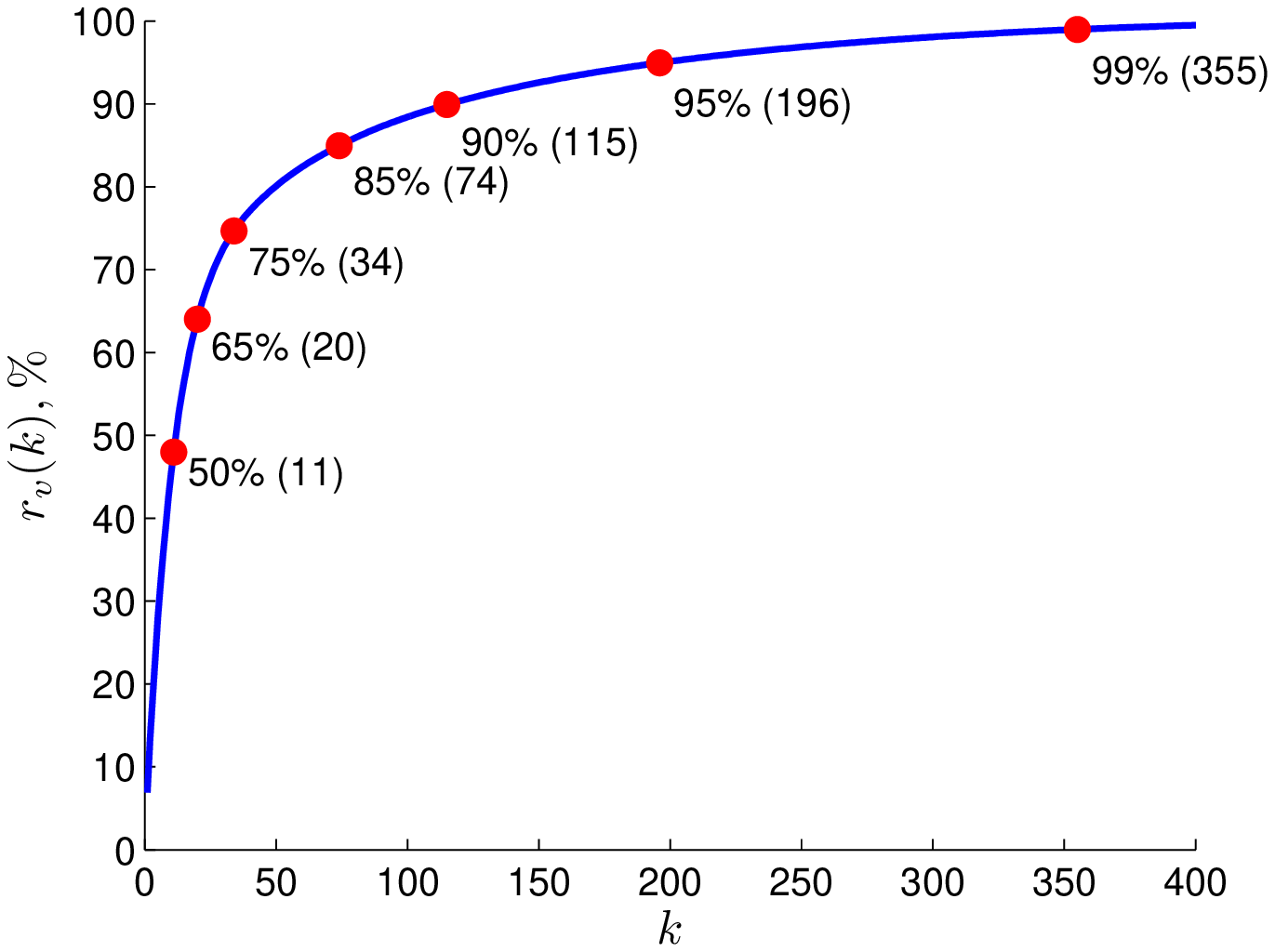}}}
  \end{center}
  \caption{(a) The values the of first 400 out of 500 eigenvalues $\lambda_k$ of matrix $Y$ obtained from a set of
    $N_r = 500$ realizations $\sigma_j$. (b)~Accumulated portion $r_v(k)$ of the variance (energy) contained in
    eigenvalues $\lambda_1$ through $\lambda_k$ as a function of $k$. For both graphs numbers in parentheses are
    the numbers $k$ of eigenvalues ($\lambda_1, \ldots, \lambda_k$) contained in the portion $r_v$ of the accumulated
    variance.}
  \label{fig:pca_eig}
\end{figure}

In this section we provide a discussion on choosing optimal number of principal components $N_{\xi}$ to reduce
dimensionality of the solution space for control $\sigma$ as discussed previously in Section~\ref{sec:PCA} in
order to optimize overall performance of our computational framework. Following this discussion, a set of
$N_r = 500$ realizations $\sigma_j$ is used to construct linear transformation matrix $\Phi^{N_{\sigma} \times N_{\xi}}$
in \eqref{eq:Phi_trunc} based on truncated SVD factorization of matrix $Y$. Figure~\ref{fig:pca_eig}(a) shows the first
400 out of 500 eigenvalues $\lambda_k$ of matrix $Y$.

Various approaches can be used to determine the size of the $\xi$-space; i.e., the $N_{\xi}$ value. Options include
the Kaiser criterion $\lambda \geq 1$ shown in \cite{Kaiser60}, the scree test introduced in \cite{Cattell66},
and the inclusion of a prescribed portion $r_v$ of the variance (energy) contained in eigenvalues $\lambda_1$
through $\lambda_k$ shown in Figure~\ref{fig:pca_eig}(b). With this last approach, given the (prescribed) parameter
$r_{\rm opt}$, $N_{\xi}$ is determined such that the following condition is satisfied
\begin{equation}
  r_v(N_{\xi}) = \dfrac{\sum_{i=1}^{N_{\xi}} \lambda_i}{\sum_{i=1}^{N_{\rm min}} \lambda_i} \cdot 100\% \geq r_{\rm opt}.
  \label{eq:pca_rv}
\end{equation}

\begin{figure}[htb!]
  \begin{center}
  \mbox{
  \subfigure[]{\includegraphics[width=0.5\textwidth]{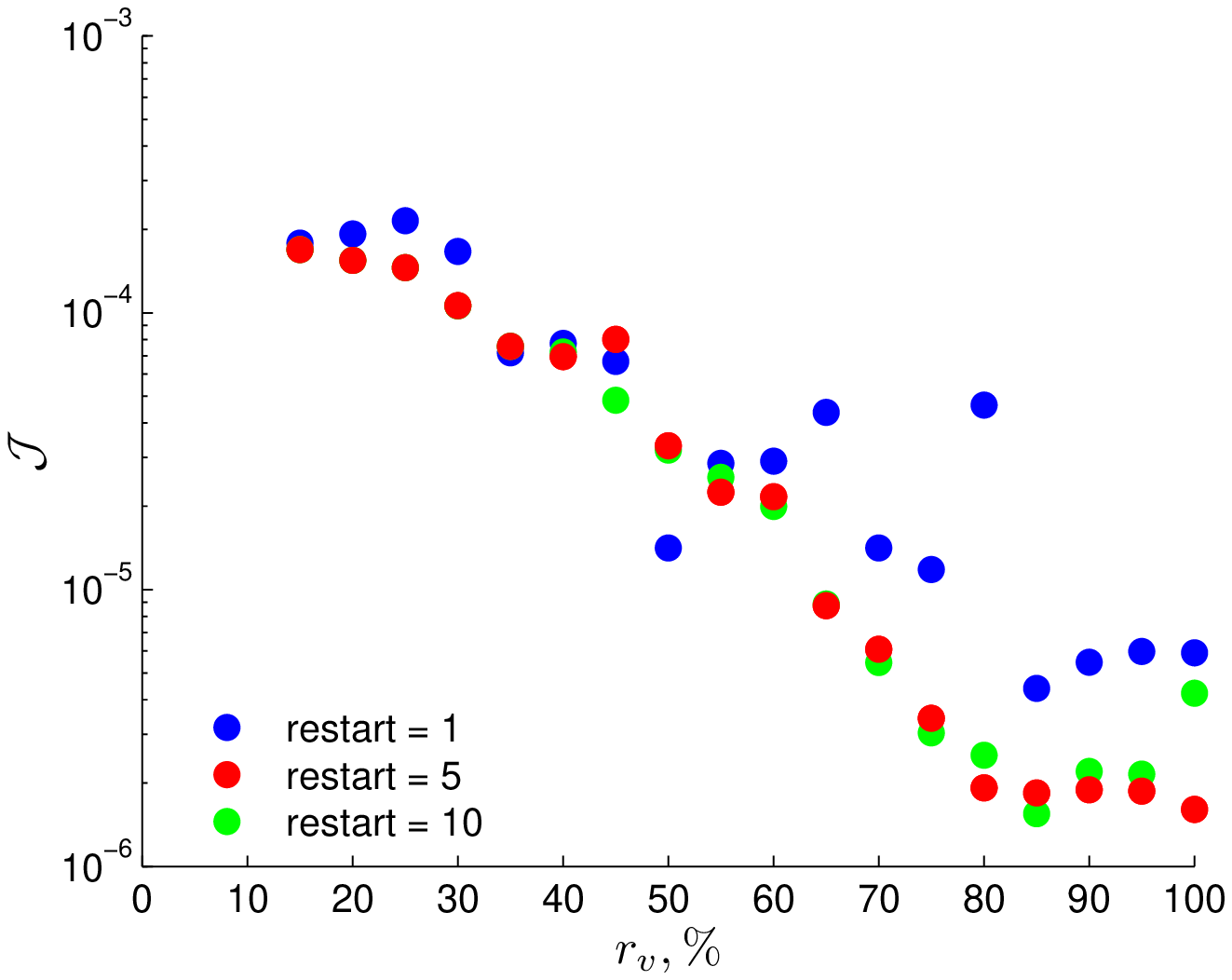}}
  \subfigure[]{\includegraphics[width=0.5\textwidth]{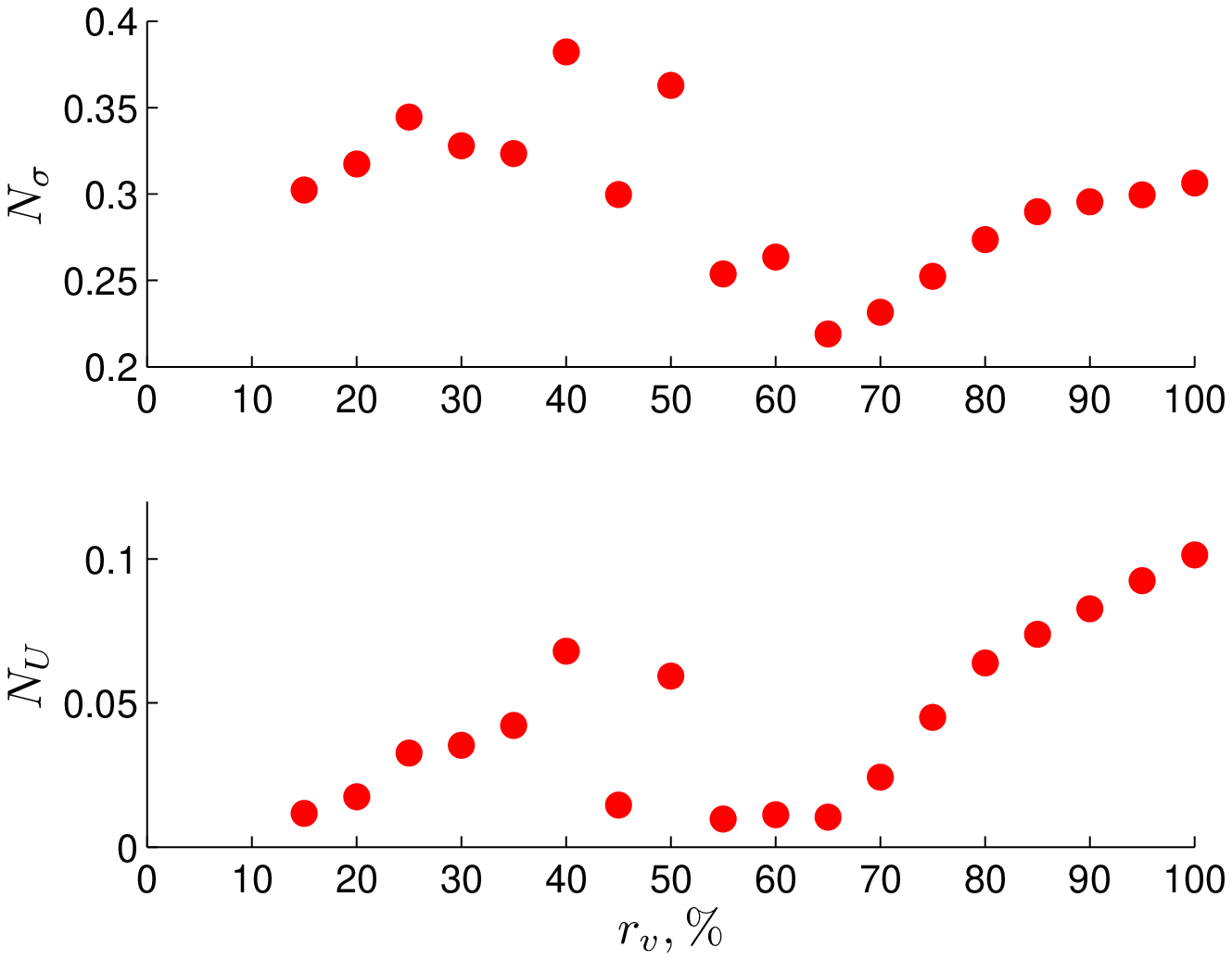}}}
  \end{center}
  \caption{(a) Cost functional $\mJ$ values evaluated at termination when Hessian approximations are restarted every
    (blue dots) $N = 1$, (red dots) $N = 5$, and (green dots) $N = 10$ iterations.
    (b) Solution norms $N_{\sigma}$ and %= \frac{\| \sigma - \sigma_{true}\|_{L_2}}{\| \sigma_{true}\|_{L_2}}$ and
    $N_U$ %= \frac{| U - U^*|}{| U^*|}$ 
    evaluated at termination for case $N=5$.}
  \label{fig:pca_opt}
\end{figure}

The scree test returns the value $r_v$ close to 80\%--85\% as at this values the graph of $r_v(k)$ is bending.
To determine parameter $r_{\rm opt}$ we run our 2D model described in Section~\ref{sec:model} 
multiple times changing the size of the $\xi$-space by setting $r_v$ value in \eqref{eq:pca_rv} to different numbers within
the range from $15\% (N_\xi = 2)$ to $100\% (N_\xi = 495)$ with step $5\%$. The performance is evaluated first by examining
cost functional $\mJ$ values at termination for three different cases to restart limited-memory quasi-Newton approximations
for Hessian in {\tt SNOPT}: every $N = 1, 5, 10$ iterations. The outcomes are represented respectively by blue, red and
green dots in Figure~\ref{fig:pca_opt}(a). The results of two cases with restarts every 5 and 10 iterations are consistent
with the scree test. We additionally examine $\sigma$ and $U$ solution
norms $N_{\sigma} = \frac{\| \sigma - \sigma_{true}\|_{L_2}}{\| \sigma_{true}\|_{L_2}}$ and
$N_U = \frac{| U - U^*|}{| U^*|}$ with results for $N = 5$ presented in Figure~\ref{fig:pca_opt}(b).
This test reveals the optimal value for $r_v$ to be close to 65\% providing only $N_{\xi} = 20$ dimensions for $\xi$-space.
This creates a high possibility for $\xi$ control space to be under-parameterized. Therefore, for all computations
shown in Section~\ref{sec:EIT_forward_backward}, unless stated otherwise, we  used $N_{\xi} = 74$ utilizing $r_{opt} = 85\%$
of accumulated variance.

\end{document}